\documentclass[11pt,twoside,leqno]{article}
\usepackage{theorem,amsfonts,amsmath,amssymb}

\newtheorem{proposition}{Proposition}[section]
\newtheorem{corollary}[proposition]{Corollary}
\newtheorem{theorem}[proposition]{Theorem}
\newtheorem{lemma}[proposition]{Lemma}
\newtheorem{remark}[proposition]{Remark}

\newtheorem{definition}[proposition]{Definition}

\newcommand{\E}{\mathbb{E}}

\newcommand{\N}{\mathbb{N}}

\newcommand{\R}{\mathbb{R}}
\newcommand{\Z}{\mathbb{Z}}
\newcommand{\Ss}{\mathbb{S}}

\newcommand{\D}{\displaystyle}

\newcommand{\proof}{\textsc{Proof: }}

\newcommand{\arcosh}{\mbox{\rm arcosh}}

\newcommand{\convrad}{\mathrm{convrad}\,}

\newcommand{\volz}{\mathrm{vol}_2}
\newcommand{\grad}{\mathrm{grad}\,}
\newcommand{\graph}{\mathrm{graph}\,}

\newcommand{\injrad}{\mathrm{injrad}\,}

\newcommand{\Int}{\mathrm{Int}}

\newcommand{\da}{\downarrow}
\newcommand{\ua}{\uparrow}

\newcommand{\ve}{\varepsilon}

\def\bi{\begin{itemize}}
\def\ei{\end{itemize}}
\def\bt{\begin{tabular}}
\def\et{\end{tabular}}

\def\kasten{\hfill\null\nobreak\hfill \hbox{\vrule\vbox{\hrule width6pt\vskip6pt\hrule}\vrule}
\par\smallskip}
\begin{document}

%\begin{flushright}
%submitted to
%\end{flushright}

\begin{center}
{\huge \bf An area bound for surfaces \vspace{0.2cm}\\ in Riemannian manifolds}
\\ \ \\
Victor Bangert\\
Mathematisches Institut der
Albert-Ludwigs-Universit\"at Freiburg,\\
Ernst-Zermelo-Stra\ss e 1, D-79104 Freiburg, Germany\\
email: victor.bangert@math.uni-freiburg.de \\
\ \\
Ernst Kuwert\\
Mathematisches Institut der
Albert-Ludwigs-Universit\"at Freiburg,\\
Ernst-Zermelo-Stra\ss e 1, D-79104 Freiburg, Germany\\
email: ernst.kuwert@math.uni-freiburg.de \\
\end{center}
\vspace{1cm}

\begin{quote}

{{\bf Abstract: } Let $M$ be a compact Riemannian manifold not containing any totally geodesic surface. Our main result shows that then the area of any complete surface immersed into $M$ is bounded by a multiple of its extrinsic curvature energy, i. e. by a multiple of the integral of the squared norm of its second fundamental form.}
\ \\ \ \\
%{\bf Keywords: } \\
%\ \\ \ \\
\end{quote}

%%%%%

\tableofcontents
\vspace{0.8cm}

\setcounter{equation}{0}

%%%%%%%%%%%%%%%%%%%%%%%%%%%%%%%%%%%%%%%%%%%%%%%%%%%%%%%%%%%%%%%
%\section{Introduction} \label{intro}
%%%%%%%%%%%%%%%%%%%%%%%%%%%%%%%%%%%%%%%%%%%%%%%%%%%%%%%%%%%%%%%
\section{Introduction}\label{sec1}
Let $M$ be a compact differentiable manifold of dimension $n \geq 3$, with a % given
smooth Riemannian metric $\bar{g}$. For a smooth immersion $f:F\to M$ of a surface $F$ into $M$ we consider % the
its extrinsic curvature energy
\begin{equation} \label{eqenergy} %% 1.1
	E(f) = \frac{1}{2} \int_F |A|^2\,d{\rm vol}_2.
\end{equation}
Here ${\rm vol}_2={\rm vol}_2^f$ is the measure associated to the %
Riemannian metric $g = f^*\bar{g}$ induced on $F$, $A$ is the normal-valued second fundamental
form of $f$, and $|A|$ denotes the euclidean norm of the tensor field $A$. Clearly, $E(f) \geq 0$ with equality if and only if $f$ is a totally geodesic immersion.
By the Gau{\ss} equation we have
\begin{equation} \label{eqgauss} %% 1.2
\frac{1}{2}|A|^2 = |A^\circ|^2 + K_g - K_{\bar{g}}^f = \frac{1}{2}|\vec{H}|^2 +  K_{\bar{g}}^f -K_g,
\end {equation}
where $K_g$ is the Gaussian curvature of $g$, and $K_{\bar{g}}^f$ denotes the sectional curvature of $\bar{g}$ % in the direction of
on the tangent planes of $f$. The second fundamental form is decomposed as $A = A^\circ + \frac{1}{2}\,g \otimes \vec{H}$,
where $A^\circ$ is trace-free and $\vec{H}={\rm{trace}}(A)$ is the mean curvature vector of $f$. If $f$ is an immersion of a closed surface $F$ into
euclidean $\E^n$, the functional $E(f)$ reduces essentially to twice the Willmore energy, 
more precisely by the Gau{\ss}-Bonnet theorem 
\begin{equation}
\label{energyeuclidean} %% 1.3
E(f) = \int_F |A^\circ|^2\,d{\rm vol_2} + 2\pi \chi(F) = \frac{1}{2} \int_F |\vec{H}|^2\,d{\rm vol}_2 - 2\pi \chi(F).
\end{equation}

In general, the area of an immersed surface is not bounded in terms of its energy. For example, consider the sequence 
$f_i$ of immersions into the flat torus $M = \E^n/\Z^n$ given by projecting $\lambda_i f$, where $f$ is a fixed immersion into $\E^n$ and 
$\lambda_i> 0$ goes to infinity. The scale invariance of $E$ implies that $E(f_i) = E(f) < \infty$ for all $i\in \N$, 
while the surface areas go to infinity. However, we prove as a main result that this behavior is rather special.
In the following, an immersion $f:F\to M$ is called complete if the Riemannian metric $f^* \bar{g}$ is complete.

\begin{theorem}[Area bound] \label{thmareabound} %% Theorem 1.1
Let $(M,\bar{g})$ be a compact Riemannian manifold of dimension $n \geq 3$. Assume that $(M, \bar{g})$ does not admit any complete,
totally geodesic surface immersions. Then there is a constant $C = C(M, \bar{g}) < \infty$ such that
${\rm vol}_2^{f}(F) \leq C \,E(f)$ for every complete, immersed surface $f: F \to M$.
\end{theorem}

For a generic metric $\bar{g}$ on $M$ there are no totally geodesic submanifolds of 
dimension $1<k<{\mathrm{dim\,}}M$ at all in $(M, \bar{g})$.
This was proved recently by Murphy \& Wilhelm \cite{MW18} for $n \geq 4$, for $n=3$ 
there is a sketch by Bryant \cite{Bry15}. See also \cite{LP22} where
Lytchak and Petrunin give a short proof valid for all $n$ in the appendix.\\
\\
In the more general context of $m$-varifolds in an $n$-dimensional 
Riemannian manifold, a related mass bound in terms of curvature 
energies was proved by A. Mondino \cite{Mon14}. Assuming that the 
bound fails, he applies a compactness argument to construct a nonzero 
$m$-varifold with vanishing generalized second fundamental form,
see Theorem 4.1 in \cite{Mon14}. However, it is not clear how 
that (a priori non-rectifiable) varifold relates to totally geodesic 
immersions, and the above generic nonexistence is not immediate 
for such a varifold.\\
\\
Our paper deals with the scale-invariant $L^2$ integral of the second 
fundamental form. The analysis would be simpler in the case of 
the $L^p$ integral with $p > 2$, employing Langer's local graph 
representations \cite{Lan85}, see also Breuning \cite{Breu12}.
In the critical case $p = 2$ there is an almost graphical 
description due to Simon \cite{Sim93}, and there are also local 
bilipschitz parametrizations by Toro \cite{Tor94}, H\'{e}lein 
\cite{Hel02}, and M\"uller \& \v{S}ver\'{a}k \cite{MS95}.
However, it remains unclear how to apply these results to
a seqeunce with mass going to infinity.\\ 
\\
Here we prove the following local version of Theorem \ref{thmareabound} 
that is an important step in the proofs of Theorem \ref{thmareabound} and
Theorem \ref{thrm13}.

\begin{theorem} \label{thrm12} Let $(M,\bar{g})$ be a compact Riemannian manifold of dimension $n \geq 3$. Assume that $(M, \bar{g})$ does not
contain any totally geodesic surface. Then for every $r>0$ there exists a constant $c(r)=c(r, M, \bar{g})>0$ such that
$E(f|B(p, r))\geq c(r)$ holds for every complete surface immersion $f:F\to M$ and every intrinsic metric ball $B(p, r)$ of radius $r$ on $F$.
In particular, if $E(f) <\infty$, then $F$ is compact.
\end{theorem}

If we relax the condition on $(M, \bar{g})$ in Theorem \ref{thmareabound} by allowing totally geodesic immersions of $S^2$,
we obtain the following slightly weaker consequence.
\begin{theorem} \label{thrm13}
Let $(M, \bar{g})$ be a compact Riemannian manifold of dimension $n\geq 3$. Assume that $(M, \bar{g})$ does not admit any
complete, totally geodesic immersions of connected surfaces other than $S^2$ or ${\R}P^2$. Then, for every constant $D$,
there exists a constant $C=C(M, \bar{g}, D)$ such that $E(f) < D$ implies ${\rm vol}_2^{f}(F) < C$ for every complete, connected, immersed
surface $f:F\to M$.
\end{theorem}

Our work is partly motivated by the joint paper \cite{KMS14} of Mondino, Schygulla and the second author.
For a compact, three-dimensional Riemannian manifold $(M, \bar{g})$, they consider the problem
of minimizing $E(f)$ in the class $[S^2,M]$ of immersions $f: S^2\to M$. They prove existence under the two assumptions:
\begin{align}
& E(f) < 4\pi \quad \mbox{ for some } f \in [S^2, M],\label{eqconditiona}\\ %% 1.4
& \mbox{For some minimizing sequence }f_i \in [S^2 ,M] \mbox{ the surface areas }\label{eqconditionb} \\ %% 1.5
& \nonumber {\rm vol}_2^{f_i}(S^2) \mbox{ remain bounded.}
\end{align}
The approach in \cite{KMS14} follows L. Simon \cite{Sim93}. In recent work by Guodong Wei \cite{Wei19} the result is reproved employing results
from \cite{MS95,Hel02,KL12,MR13} on conformal parametrizations. Combining the result from \cite{KMS14} with Theorem \ref{thrm13} we obtain
\begin{corollary} \label{coroll14}
Let $(M,\bar{g})$ be a three-dimensional, compact Riemannian manifold that admits no complete, totally geodesic immersions
of connected surfaces other than $S^2$ or ${\R}P^2$. If (\ref{eqconditiona}) is satisfied, there exists $f\in [S^2, M]$ minimizing $E$ on $[S^2,M]$.
\end{corollary}

Condition (\ref{eqconditiona}) holds if $Scal_{\bar{g}}(x) >0$ for some $x \in M$. Condition (\ref{eqconditionb}) is an easy consequence
of the Gau{\ss} equations when $(M,\bar{g})$ has positive sectional curvature. In particular, \cite{KMS14} proves the existence of a minimizer in
$[S^2, M]$, if $M$ is compact and has positive sectional curvature. We recover this result, since, by the Bonnet-Myers and the Gau\ss-Bonnet theorems,
a complete, connected, totally geodesic immersed surface in a compact manifold of positive sectional curvature is of type $S^2$ or ${\R}P^2$.\\

Next we shortly outline the proofs of Theorems \ref{thmareabound}--\ref{thrm13}. These proofs are by contradiction.
To prove Theorem \ref{thrm12} we assume that there exists a sequence $f_i: F_i \to M$ of complete surface immersions and a sequence of intrinsic metric balls
$B(p_i, r)\subseteq F_i$ of fixed radius $r>0$ such that $\lim_{i\to \infty}E(f_i|B(p_i, r))=0$. Then we prove that the set of limit points of sequences $f_i(x_i)$
with $x_i\in B(p_i, r)$ contains a totally geodesic surface, see Theorem \ref{thrm94}, and Theorem \ref{thrm131} for a statement concerning Hausdorff convergence.
This proof relies on the results of Sections \ref{sec2}--\ref{sec8} that will be reviewed below. We will also consider the following stronger assumption:
\begin{equation} \label{eq16}
\mbox{There exist sequences } p_i\in F_i, R_i\to \infty, \mbox{such that } \lim_{i\to \infty}E(f_i|B(p_i, R_i))=0.
\end{equation}

Under assumption (\ref{eq16}) we can prove that $M$ admits a  {\em complete}, totally geodesic surface immersion, see Corollary \ref{coroll105}.
To derive Theorem \ref{thmareabound} from this last result we use a Voronoi type covering argument to verify assumption (\ref{eq16}). Here we need to
generalize the well-known Bol-Fiala type area estimates from parallel sets to weakly starshaped sets, cf. Section \ref{sec11}. In Section \ref{sec13} we prove
Theorem \ref{thrm13} by contradiction. Here assumption (\ref{eq16}) is easily seen to be satisfied, so that we obtain a complete, totally geodesic surface immersion
$f: N\to M$. Then we show that $N$ is not diffeomorphic to $S^2$ or $\R P^2$. This is inspired by G.\,Reeb's stability theorem \cite{Reeb47} % [25]
from the theory of foliations.\\

Now we review the contents of Sections \ref{sec2}--\ref{sec8} on which the proofs of Theorems \ref{thmareabound}--\ref{thrm13} are based. In Section \ref{sec2} we use
Gronwall's inequality to conclude that arclength-parametrized curves $\gamma:[a, b]\to M$ with small total absolute geodesic curvature are $C^1$-close to geodesics.
In Section \ref{sec3} we consider a complete immersion $f:F\to M$ such that $E(f)/{\rm vol}_k^{f}(F)$ is small where $k={\rm dim} F$. We use the invariance of the
Liouville measure under the geodesic flow to find a subset of the unit tangent bundle $SF$ of large Liouville measure such that the $f$-images of geodesics in $F$
with initial vectors in this subset are $C^1$-close to geodesics in $M$. More generally, we prove a similar statement for configurations $(c_1, c_2, t)$ where
$c_1, c_2$ are geodesics in $F$ with $c_2(0) = c_1(t)$. To relate the estimates on the Liouville measure to the geometry of $F$, we need lower estimates on the volume of intrinsic metric balls in $F$, see Sections \ref{sec6} and \ref{sec8}. While the results in Sections \ref{sec2} and \ref{sec3} are true for manifolds $F$ of arbitrary dimension we can prove these volume estimates and much of the following only if ${\rm dim}F=2$. This involves some new results on the geometry of surfaces in euclidean space $\E^n$ that depend on bounds on their energy, see Sections \ref{sec4} and \ref{sec5}. In Section \ref{sec7} we consider a sequence of complete surface immersions $f_i:F_i\to M$, and assume the existence of a sequence of balls $B(p_i, R_i)\subseteq F_i$ such that
$R_i\to \infty$ and $\liminf_{i\to\infty} E(f_i|B(p_i, R_i))<\frac{\pi}{4}$. Relying on ideas by T. Shioya \cite{Shio99} we prove that a subsequence of the sequence
$(F_i, p_i)$ converges to a proper, pointed length space $(Y, y_0)$ with respect to pointed Gromov-Hausdorff convergence and that $Y$ has locally finite 2-dimensional Hausdorff measure. Moreover, we can assume that the immersions $f_i$ converge to a 1-Lipschitz map $f:Y\to M$. On the other hand, the results of Section \ref{sec3} can be used to see that the Hausdorff dimension of $f(Y)$ is at least three, unless $f(Y)\subseteq M$ contains a totally geodesic surface, see Section \ref{sec9}.\\

{\bf Notation.} Here we collect some of the notation used throughout the paper. The unit sphere of a euclidean vector space $\E$ is denoted by $S\E$.
If ${\rm dim}\E=k$ then $\alpha_{k-1}$ is the $(k-1)$-volume of $S\E$. If $M$ is a manifold and $2\leq k\leq{\rm dim}M$ then $\pi_G: G_k M\to M$ denotes the Grassmann bundle of $k$-dimensional linear subspaces in $TM$. For Riemannian manifolds $(M, g)$, we let $\pi: SM \to M$ denote the unit sphere bundle. The intrinsic metric ball in $M$ with center $p\in M$ and radius $r>0$ is denoted by $B(p, r)$. If ${\rm dim}M = m$ then ${\rm vol}_m$ is the Riemannian volume on $M$, and ${\cal H}^k$, $1\leq k\leq {\rm dim }M$, denotes $k$-dimensional Hausdorff measure on $M$. Geodesics will always be parametrized by arclength, and $c_v$ will denote the geodesic with initial vector $v\in SM$. If $\gamma$ is a curve in $M$ then ${\cal P}_t^\gamma: T_{\gamma(0)} M \to T_{\gamma(t)} M$ denotes parallel translation along
$\gamma$ from $\gamma(0)$ to $\gamma(t)$. The injectivity radius of $M$ will be denoted by ${\rm{injrad}\,}(M)$.\\
\\
{\bf Acknowledgements: }We thank C. Debin for communication concerning \cite{Deb20}. We thank
the referee for his/her careful reading.

\setcounter{equation}{0}
%%%%%%%%%%%%%%%%%%%%%%%%%%%%%%%%%%%%%%%%%%%%%%%%%%%%%%%%%%%%%%%
\section{$L^1$-almost geodesics} \label{sec2}
Let $(M, g)$ denote a compact Riemannian manifold with Levi-Civit\`{a} connection $\nabla$. In this section we consider arclength-parametrized curves
in $M$ for which the $L^1$-norm of the covariant derivative of $\dot{\gamma}$ is small. Using Gronwall's inequality we will easily see that such
curves are $C^1$-close to geodesics. This will be applied to geodesics of submanifolds of $M$ along which the norm of the second fundamental form has small integral.\\

On the unit tangent bundle $SM$ we consider the Sasaki metric $\bar{g}$ induced by $g$, see \cite{Sas58}, and the distance function $d^{SM}$ induced by $\bar{g}$.

\begin{lemma}\label{lemma21}
There exist constants $B>0$, $C>0$, such that the following holds for all arclength-parametrized $C^2$-curves
$\gamma:[a, b]\to M$. If $c:[a,b]\to M$ is the geodesic with initial vector $\dot{c}(a) = \dot{\gamma}(a)$ and $t\in [a, b]$, then
$$
d^{SM}(\dot{\gamma}(t), \dot{c}(t)) \leq B e^{C(t-a)} \int_{a}^{t}|\nabla_{\frac{\partial}{\partial t}} \dot{\gamma}| (s)\, ds.
$$
\end{lemma}

\proof We recall the following qualitative version of Gronwall's inequality. Let $(N, h)$ be a compact Riemannian manifold with induced distance $d^N$,
and let $X$ be a $C^1$-vector field on $N$ with flow $\Phi:N\times \R\to N$. Then there exist constants $B>0$, $C>0$ such that the following holds
for all $C^1$-curves $\beta:[a,b]\to N$ and all $t\in[a,b]$:
\begin{equation}\label{eq21}
d^{N}\big(\beta(t), \Phi_t(\beta(a))\big) \leq B e^{C(t-a)} \int_{a}^{t}|\dot{\beta}(\tau)- X(\beta(\tau))|\, d\tau.
\end{equation}

We apply this to the case $N=SM$, $h=\bar{g}$, and to the vector field $X$ on $SM$ whose flow $\Phi$ is the geodesic flow on $SM$.
Given an arclength-parametrized $C^2$-curve $\gamma:[a, b]\to M$, we consider $\beta:[a, b]\to SM$, $\beta(t)=\dot{\gamma}(t)$.
By the definition of the Killing metric $\bar{g}$, we have
$$
|\dot{\beta}(t) - X(\beta(t))|^{\bar{g}}=|\nabla_{\frac{\partial}{\partial t}} \dot{\gamma}(t)|^g.
$$
Hence \ref{eq21} implies our claim. \kasten\vspace{0.2cm}

The following statement is a simple consequence of the dependence of solutions of linear ODEs on the equation.

\begin{lemma} \label{lemma22}
Suppose $\gamma_i: [-R, R]\to M$ is a sequence of $C^1$-curves that converge to $\gamma: [-R, R]\to M$ uniformly in the $C^1$-topology.
Let $W_i$ be parallel vector fields along $\gamma_i$ such that $\lim_{i\to \infty} W_i(0) = w\in T_{\gamma(0)} M$.
Then the $W_i$ converge uniformly to the parallel vector field $W$ along $\gamma$ with $W(0) = w$.
\end{lemma}

These lemmas will be applied in the following situation. We consider a sequence of complete immersions $f_i: F_i\to M$ of $k$-dimensional manifolds $F_i$ into $M$.
The second fundamental forms of the $f_i$ will be denoted by $A_i$.

\begin{proposition} \label{prop23}
Let $v_i\in SF_i$ be a sequence such that $\lim_{i\to \infty} df_i(v_i) = v\in SM$ exists. Assume that $R>0$ and
$\lim_{i\to \infty} \int_{-R}^R |A_i|\circ c_{v_i}(t)\,dt = 0$. Then the following statements {\em (a)--(c)} are true.
\begin{itemize}
\item[\em{(a)}] The curves $f_i\circ c_{v_i}|[-R, R]$ converge to the geodesic $c_v|[-R, R]$ uniformly in the $C^1$-topology.
\item[\em{(b)}] If $W_i:[-R, R]\to TF_i$ are parallel vector fields along $c_{v_i}|[-R, R]$ and $\lim_{i\to \infty} df_i(W_i(0)) = w\in TM$ exists,
           then the vector fields $df_i \circ W_i$ along $f_i\circ c_{v_i}|[-R, R]$ converge uniformly to the parallel vector field $W$ along $c_v|[-R, R]$ with $W(0) = w$.
\item[\em{(c)}] If $\lim_{i\to \infty} df_i(T_{c_{v_i}(0)}F_i) = \bar{L}\in G_k M$ exists, then the curves $L_i: [-R, R]\to G_k M$, $L_i(t)= df_i(T_{c_{v_i}(t)}F_i)$,
           converge uniformly to the curve $L: [-R, R]\to G_k M$, where $L(t) = {\cal{P}}_t^{c_v}(\bar{L})$ is the parallel translate of $\bar{L}$ along $c_v$.
\end{itemize}
\end{proposition}

\proof
\begin{itemize}
\item[(a)] By the definition of the second fundamental form $A_i$, we have
$$
\nabla_{\frac{\partial}{\partial t}} (f_i\circ c_{v_i})^{\dot{}}= A_i(\dot{c}_{v_i}, \dot{c}_{v_i}).
$$
Hence (a) is a consequence of Lemma \ref{lemma21} and the convergence of the $(M, g)$-geodesics $c_{df_i(v_i)}|[-R, R]$ to $c_v|[-R, R]$.

\item[(b)] Let $Z_i:[-R, R] \to TM$ denote the parallel vector fields along $f_i\circ c_{v_i}|[-R, R]$ with $Z_i(0)=df_i(W_i(0))$. Then we have
$$
\langle df_i \circ W_i, Z_i\rangle (t) = \langle df_i \circ W_i, Z_i\rangle (0)+ \int_0^t \langle A_i (\dot{c}_{v_i}, W_i), Z_i\rangle (\tau)\,d\tau.
$$
Hence our assumption $\lim_{i\to \infty} \int_{-R}^R |A_i|\circ c_{v_i}(t)\,dt = 0$ implies that $\langle df_i \circ W_i, Z_i\rangle$ converges uniformly on $[-R, R]$
to the constant $|Z_i(0)|^2 = |W_i(0)|^2$. Since $|(df_i\circ W_i)(t)|=|W_i(0)|=|Z_i(t)|$, this implies $\lim_{i\to \infty} |df_i\circ W_i - Z_i|=0$
uniformly on $[-R, R]$. Using (a), our assumption $\lim_{i\to \infty} df_i(W_i(0))= w = W(0)$, and Lemma \ref{lemma22}, we see that the $Z_i$, and hence
$df_i\circ W_i$, converge uniformly to $W$.
\item[(c)] This is a direct consequence of (b).
\end{itemize}
\vspace{-0.2cm} \kasten

\begin{remark} \label{rem24}
{\em In the applications of Proposition \ref{prop23} in Sections \ref{sec7}, \ref{sec9}, \ref{sec10} and \ref{sec13}, the hypothesis
$\lim_{i\to \infty} \int_{-R}^R |A_i|^2\circ c_{v_i}(t)\,dt = 0$ will hold. By the Cauchy-Schwarz inequality this implies
$\lim_{i\to \infty} \int_{-R}^R |A_i|\circ c_{v_i}(t)\,dt = 0$.}
\end{remark}

\setcounter{equation}{0}
%%%%%%%%%%%%%%%%%%%%%%%%%%%%%%%%%%%%%%%%%%%%%%%%%%%%%%%%%%%%%%%
\section{Estimates arising from integral geometry} \label{sec3}

We consider a $k$-dimensional, complete Riemannian manifold $(F, g)$. The unit tangent bundle of $F$ will be denoted by $\pi: SF\to F$
with fibers $S_p F=\pi^{-1}(p)$ for $p\in F$. Our estimates will follow from the invariance of the Liouville measure ${\cal{L}}$ on $SF$
under the geodesic flow $\Phi: SF\times \R\to SF$. We recall that locally the Liouville measure ${\cal{L}}$ is the product of the Riemannian volume $\mathrm{vol}_k$
on $F$ and the standard $(k-1)$-volume $\mathrm{vol}_{S_p F}$ on the euclidean spheres $S_p F$. For $v\in SF$, we let $c_v: \R\to F$
denote the geodesic with initial vector $\dot{c}_v(0) = v$, i.\,e. $c_v(t)=\pi\circ\Phi(v, t)$. ``Measurability'' will be understood with respect to the Borel $\sigma$-algebra.
Lebesgue measure on $\R$ will be denoted by $\lambda$.\\

\begin{lemma}\label{lemma31}
Let $h: F\to [0, \infty]$ be measurable. For arbitrary $\ve >0$, $R>0$ consider the set
$$
{\cal{V}}_{\ve, R}(h)= \Big\{v\in SF| \int_{-R}^R h\big(c_v(t)\big)\,dt \geq \ve\Big\}.
$$
Then we have
$$
{\cal{L}}\big({\cal{V}}_{\ve, R}(h)\big) \leq \alpha_{k-1} \frac{2R}{\ve} \int_F h\, d\mathrm{vol}_k,
$$
where $\alpha_{k-1}$ denotes the $(k-1)$-volume of the unit sphere in euclidean $k$-space.
\end{lemma}

\proof
The invariance of ${\cal{L}}$ under $\Phi$ implies that the integral
$$
\int_{SF} h\big(c_v(t)\big)\, d{\cal{L}}(v)
$$
is independent of $t\in \R$. From the definition of ${\cal{L}}$ we obtain
$$
\int_{SF} h\big(c_v(t)\big)\, d{\cal{L}}(v) = \int_{SF} h\big(c_v(0)\big)\, d{\cal{L}}(v) =
\alpha_{k-1} \int_F h\,d\mathrm{vol}_k.
$$
This implies
$$
{\cal{L}}\big({\cal{V}}_{\ve, R}(h)\big) \cdot \ve \leq
\int_{SF} \bigg(\int_{-R}^R h\big(c_v(t)\big)\, dt\bigg) d{\cal{L}}(v) =
\alpha_{k-1} 2R \int_F h\, d\mathrm{vol}_k.
$$
\kasten\vspace{0.2cm}

\begin{lemma}\label{lemma32}
If $B\subseteq F$ is a Borel set and $R>0$, then
$$
({\cal{L}}\times \lambda) (\{(v, t)\in SF\times [-R, R]|c_v(t)\in B\}=\alpha_{k-1}2R \mathrm{vol}_k(B).
$$
\end{lemma}

\proof As in the preceding proof we see that ${\cal{L}}(\{v\in SF |c_v(t)\in B\})$ is independent of $t\in\R$.
\kasten\vspace{0.2cm}

We will apply Lemma \ref{lemma31} in the following situation.
We will consider an isometric immersion $f: F\to M$ of $F$ into a Riemannian manifold $M$,
and let $h$ be the squared norm $|A|^2$ of the second fundamental form of $f$.
Under appropriate conditions on $E(f)=\frac{1}{2} \int_F|A|^2 d\mathrm{vol}_k$, we can use Lemma \ref{lemma31} to find a large set of vectors
$v\in SF \backslash {\cal{V}}_{\ve, R}(|A|^2)$, i.\,e. vectors $v\in SF$ for which $\int_{-R}^R |A|^2\circ c_v(t)\,dt <\ve$.
If additionally $\sqrt{2R\ve}$ is small, then Proposition \ref{prop23} shows that for these vectors $v\in SF$ the curve
$f\circ c_v |[-R, R]$ is $C^1$-close to a geodesic in $M$.\\

In the proof of the existence of totally geodesic surfaces in Section \ref{sec9} we will need pairs of geodesics $c, \tilde{c}$
in $F$ such that $\tilde{c}(0) = c(t)$ for some $t\in\R$ and such that both $f\circ c$ and $f\circ \tilde{c}$ are $C^1$-close to geodesics in $M$ on appropriate intervals $[-R, R]$ resp. $[-r, r]$. We will encode such pairs of geodesics as follows. We consider $SF^2 = \bigcup_{p\in F}(S_pF\times S_pF)$ and the bundle
$\tilde{\pi}: SF^2\times \R\to F$, $\tilde{\pi}(v, w, t)=\pi(v)=\pi(w)$. ${\cal{L}}^2$ will denote the measure on $SF^2$ that is the product of $\mathrm{vol}_k$ with two factors $\mathrm{vol}_{S_pF}$. A tuple $(v, w, t)\in SF^2\times \R$ encodes the pair of geodesics $c=c_v$ and $\tilde{c}$ with $\dot{\tilde{c}}(0)={\cal{P}}_t^{c_v}(w)$, where
${\cal{P}}_t^{c_v}: S_{c_v(0)}F\to S_{c_v(t)}F$ denotes parallel transport along $c_v$.\\

Given $\ve>0$, $R>0$, $r>0$, we set
\begin{equation} \label{eq31}
{\cal{V}}_{\ve,R,r}^2(h) =\{(v, w, t)\in SF^2\times [-R, R]| v\in {\cal{V}}_{\ve, R}(h) \mbox{ or }
{\cal{P}}_t^{c_v}(w)\in {\cal{V}}_{\ve, r}(h)\}.
\end{equation}

So $(v, w, t)\in (SF^2\times [-R, R])\backslash {\cal{V}}_{\ve,R,r}^2(h)$ implies $\int_{-R}^R |A|^2 \circ c_v(t)\,dt < \ve$ and
$\int_{-r}^r |A|^2 \circ c_{\tilde{w}}(t)\,dt < \ve$, where $\tilde{w}={\cal{P}}_t^{c_v}(w)$.

\begin{lemma}\label{lemma33}
$({\cal{L}}^2\times \lambda)({\cal{V}}_{\ve,R,r}^2(h)) \leq \alpha_{k-1}^2 \frac{4R}{\ve} (R+r) \int_F h\,d\mathrm{vol}_k$.
\end{lemma}

\proof We will prove that
\begin{equation} \label{eq32}
({\cal{L}}^2\times \lambda)\{(v, w, t)\in SF^2\times [-R, R]| v\in {\cal{V}}_{\ve, R}(h)\}=\alpha_{k-1} 2R\,{\cal{L}}\big({\cal{V}}_{\ve,R}(h)\big)
\end{equation}
and that
\begin{equation} \label{eq33}
({\cal{L}}^2\times \lambda)\{(v, w, t)\in SF^2\times [-R, R]| P_t^{c_v}(w)\in {\cal{V}}_{\ve, r}(h)\}=\alpha_{k-1} 2R {\cal{L}}\big({\cal{V}}_{\ve,r}(h)\big).
\end{equation}
Then our claim will follow by combining (\ref{eq32}), (\ref{eq33}) with Lemma \ref{lemma31}.
Equation (\ref{eq32}) is a direct consequence of Fubini's theorem. To prove equation (\ref{eq33}) note that the orthogonality of
${\cal{P}}_t^{c_v}: T_{\pi(v)}F \to T_{c_v(t)}F$ and the invariance of ${\cal{L}}$ under the geodesic flow imply that, for every $t\in \R$, we have
$$
\begin{array}{c}
\D{\cal{L}}^2(\{(v, w)\in SF^2|{\cal{P}}_t^{c_v}(w)\in {\cal{V}}_{\ve, r}(h)\})=\int_{SF}\mathrm{vol}_{S_{c_v(t)} F}({\cal{V}}_{\ve,r}(h)\cap S_{c_v(t)})\,d{\cal{L}}(v)\\
\D=\int_{SF}\mathrm{vol}_{S_{\pi(v)}F}({\cal{V}}_{\ve,r}(h)\cap S_{\pi(v)}F)\,d{\cal{L}}(v)=\alpha_{k-1}{\cal{L}}\big({\cal{V}}_{\ve,r}(h)\big).
\end{array}
$$
Integrating this equation over $t\in[-R, R]$ we obtain (\ref{eq33}).
\kasten\vspace{0.2cm}

We will also need the following localized version of Lemma \ref{lemma33}.
\begin{lemma}\label{lemma34}
$({\cal{L}}^2\times \lambda)\big({\cal{V}}_{\ve,R,r}^2(h)\cap \tilde{\pi}^{-1} (B(p, r))\big)\leq
\alpha_{k-1}^2 \frac{4R}{\ve}(R+r)\int_{B(p, R+2r)} h\,d\mathrm{vol}_k$.
\end{lemma}

\proof We consider $\tilde{h}= h\cdot \chi_{B(p,R+2r)}$, and show that
$\big({\cal{V}}_{\ve,R,r}^2(h)\cap \tilde{\pi}^{-1} (B(p, r))\big)\subseteq {\cal{V}}_{\ve,R,r}^2(\tilde{h})$.
Indeed, if $(v, w, t)\in {\cal{V}}_{\ve,R,r}^2(h)\cap \tilde{\pi}^{-1} (B(p, r))$ and $v\in {\cal{V}}_{\ve,R}(h)$, then
$c_v(t)\in {B(p,R+r)}$ for all $|t|\leq R$, hence $h\circ c_v(t)=\tilde{h}\circ c_v(t)$ for $|t|\leq R$, so that $v\in {\cal{V}}_{\ve,R}(\tilde{h})$.
If $(v, w, t)\in {\cal{V}}_{\ve,R,r}^2(h)\cap \tilde{\pi}^{-1} (B(p, r))$ and ${\cal{P}}_t^{c_v}(w)\in {\cal{V}}_{\ve, r}(h)$, then $c_v(t)\in {B(p,R+r)}$
and, as above, we see that ${\cal{P}}_t^{c_v}(w)\in {\cal{V}}_{\ve, r}(\tilde{h})$.
\kasten\vspace{0.2cm}

We define the measurable function $H= H_{\ve, R, r}(h):F\to [0, \alpha_{k-1}^2 2R]$ as the fibrewise 
volume of ${\cal{V}}_{\ve,R,r}^2(h)$, i.\,e.
\begin{equation} \label{eq34}
H(p)= (\mathrm{vol}_{S_p F}\times \mathrm{vol}_{S_p F}\times \lambda) 
\big({\cal{V}}_{\ve, R, r}^2(h)\cap (S_p F\times S_p F\times\R)\big).
\end{equation}

Then Lemma \ref{lemma33} resp. Lemma \ref{lemma34} imply
\begin{equation} \label{eq35}
\int_F H\, d\mathrm{vol}_k \leq \alpha_{k-1}^2 \frac{4R}{\ve} (R+r) \int_F h\,d\mathrm{vol}_k,
\end{equation}
and
\begin{equation} \label{eq36}
\int_{B(p, r)} H\,d\mathrm{vol}_k\leq \alpha_{k-1}^2 \frac{4R}{\ve} (R+r) \int_{B(p, R+2r)} h\,d\mathrm{vol}_k.
\end{equation}

\begin{definition} \label{def35}
The set $G_{\ve, R,r}(h)$ of ``$(h, \ve, R, r)$-good points'' consists of all $p\in F$ such that
$H(p) < \ve$, i.\,e. $G_{\ve, R, r}(h)= H^{-1}([0, \ve))$.
\end{definition}
Note that $p\in G_{\ve, R,r}(h)$ iff -- up to a set of measure smaller than $\ve$ -- the tuples $(v, w, t)\in S_pF \times S_pF \times [-R, R]$
encode geodesics $c=c_v$ and $\tilde{c}$, $\dot{\tilde{c}}(0) = {\cal{P}}_t^{c_v}(w)$, such that
$\int_{-R}^R h\circ c(t)\,dt < \ve$ and $\int_{-r}^r h \circ \tilde{c}(t)\,dt < \ve$. In particular, if $r\leq{\tilde{r}}$ and $R\leq{\tilde{R}}$, then
$G_{\ve, \tilde{R}, \tilde{r}}(h)\subseteq G_{\ve, R,r}(h)$. This obvious inclusion will be used at various instances without further notice.\\

The following proposition is a consequence of inequalities (\ref{eq35}) and (\ref{eq36}).
The following estimates hold for all $p\in F$, all $\ve>0$, $R>0$, $r>0$, and all measurable functions $h:F\to[0, \infty]$.

\begin{proposition}\label{prop36}
{\em{(a)}} $\D\mathrm{vol}_k (F\backslash G_{\ve, R, r}(h)) \leq \alpha_{k-1}^2\frac{4R}{\ve^2} (R+r) \int_F h\,d\mathrm{vol}_k$ \\
{\em{(b)}} $\D\mathrm{vol}_k (B(p, r) \backslash G_{\ve, R, r}(h)) \leq \alpha_{k-1}^2 \frac{4R}{\ve^2} (R+r) \int_{B(p, R+2r)} h\,d\mathrm{vol}_k$.
\end{proposition}

\proof
Since $F\backslash G_{\ve, R, r}(h)) = H^{-1}([\ve, \infty))$ we have $\ve\,\mathrm{vol}_k (F\backslash G_{\ve, R, r}(h))\leq \int_F H\,d\mathrm{vol}_k$.
Hence (a) is a consequence of (\ref{eq35}). Similarly (b) follows from (\ref{eq36}).
\kasten\vspace{0.2cm}

Assuming that $\int_F h\,d\mathrm{vol}_k$ is small we want to use Proposition \ref{prop36} to conclude that $G_{\ve, R, r}(h)$ is almost dense in $F$. For this
we need a lower bound on the volume of metric balls. For $k=2$ and under appropriate assumptions on $(F, g)$, such a lower bound will be proved in Section \ref{sec6},
see Proposition \ref{prop63}.\\

While Proposition \ref{prop36} suffices for the proof of the existence of (pieces of) totally geodesic surfaces, a slightly different estimate will be useful in the
proof of the existence of complete totally geodesic surfaces in Section \ref{sec10}. Here we need not only ``good'' points $p$, but points such that additionally
$c_v(t)$ is ``good'' for most $(v, t)\in S_pF\times [-R, R]$. In this application it will not be necessary to discriminate between the roles of $R$ and $r$.
So we will set $R=r$, and abbreviate
\begin{equation} \label{eq37}
G_{\ve, R, R}(h))=G_{\ve, R}(h).
\end{equation}

For fixed $\ve>0$, $R>0$, and $h$, we define $l:SF\to [0, 2R]$ by
$$
l(v)=\lambda(\{t\in[-R, R]|c_v(t)\notin G_{\ve, R}(h)\}) = 2R-\lambda(\{t\in[-R, R]|c_v(t)\in G_{\ve, R}(h)\}).
$$
Finally we set
\begin{equation} \label{eq38}
\tilde{G}_{\ve, R}(h) =\{q\in G_{\ve, R}(h)|\,\mathrm{vol}_{S_qF}\big(l^{-1}([\ve, 2R])\cap S_qF\big)<\ve\}.
\end{equation}

Note that $q\in G_{\ve, R}(h)$ is in $\tilde{G}_{\ve, R}(h)$ iff -- up to a set of $\mathrm{vol}_{S_qF}$-measure smaller than $\ve$ -- the vectors $v\in S_qF$ satisfy
\begin{equation} \label{eq39}
\lambda\big(\{t\in[-R, R] | c_v(t)\in G_{\ve, R}(h)\}\big)\geq 2R-\ve.
\end{equation}

As a consequence of Proposition \ref{prop36} and Lemma \ref{lemma32} we obtain:
\begin{corollary} \label{coroll37}
$\D\mathrm{vol}_k (B(p, R) \backslash \tilde{G}_{\ve, R}(h)) \leq c_k(R, \ve)\int_{B(p, 5R)} h\,d\mathrm{vol}_k$,\\
where $c_k(R, \ve)= \alpha_{k-1}^2 \frac{8R^2}{\ve^2} \left(1+ \alpha_{k-1}\frac{6R}{\ve^2}\right)$.
\end{corollary}

\proof We are going to prove that
\begin{equation} \label{eq310}
\mathrm{vol}_k (B(p, R) \cap(G_{\ve, R}\backslash \tilde{G}_{\ve, R})) \leq
\alpha_{k-1}^3 \frac{48R^3}{\ve^4} \int_{B(p, 5R)} h\,d\mathrm{vol}_k.
\end{equation}
Combined with Proposition \ref{prop36}(b) this estimate implies our claim.

To prove (\ref{eq310}) we first note that
$$
\{(v,t)\in \pi^{-1}(B(p, R))\times[-R, R]|c_v(t)\notin G_{\ve, R}(h)\} \subseteq \{(v,t)\in SF\times[-R, R]|c_v(t)\in B(p, 2R)\backslash G_{\ve, R}(h)\}.
$$
Applying Lemma \ref{lemma32} to the right hand side of this inclusion and using the definitions of $l$ and
$\tilde{G}_{\ve, R}(h))$, we obtain
$$
\begin{array}{rcl}
\alpha_{k-1} 2R\,\mathrm{vol}_k \big(B(p, 2R)\backslash G_{\ve, R}(h)\big) &\geq&
({\cal{L}}\times \lambda)\big(\{(v,t)\in \pi^{-1}(B(p, R))\times[-R, R]|c_v(t)\notin G_{\ve, R}(h)\}\big)=\\
\D\int_{\pi^{-1}(B(p, R))} l(v)\, d{\cal{L}}(v) &\geq& \ve^2  \mathrm{vol}_k (B(p, R) \cap\big(G_{\ve, R}(h)\backslash \tilde{G}_{\ve, R}(h))\big).
\end{array}
$$
Now we apply Proposition \ref{prop36}(b) to the first term in this chain of inequalities and obtain (\ref{eq310}).
\kasten

%%%%%%%%%%%%%%%%%%%%%%%%%%%%%%%%%%%%%%%%%%%%%%%%%%
\setcounter{equation}{0}
%%%%%%%%%%%%%%%%%%%%%%%%%%%%%%%%%%%%%%%%%%%%%%%%%%%%%%%%%%%%%%%
\section{A lower bound for the total absolute geodesic curvature of simple closed curves on surfaces in euclidean spaces}\label{sec4}

From now on $F$ will denote a $2$-dimensional manifold, in contrast to 
the preceding section where $\mathrm{dim}\,F = k$ was arbitrary. In this 
section we consider a complete, connected surface $F$ immersed into 
euclidean space $\E^n$ and a smooth simple closed curve $\Gamma$ on $F$.
The total absolute geodesic curvature $|{\cal{K}}|(\Gamma)$ of $\Gamma$ is defined by $|{\cal{K}}|(\Gamma)= \int_{\Gamma}|\kappa_g(q)|\, ds(q)$,
where $|\kappa_g(q)|$ denotes the absolute value of the geodesic curvature of $\Gamma$ at $q\in \Gamma$.
The second fundamental form of the immersion of $F$ into $\E^n$ will be denoted by $A$. The tubular neighborhood $\Gamma^t$ of $\Gamma$ of radius $t>0$
is the set of points $x\in F$ that can be joined to $\Gamma$ by a curve on $F$ of length at most $t$, i.\,e. $\Gamma^t= \{x\in F|d^F(x, \Gamma)\leq t\}$.

\begin{proposition} \label{prop41}
Given $c\in (0, \frac{4}{3}\pi)$ there exists $\delta=\delta(c)>0$ and $\beta = \beta(c)\geq 1$ such that the following holds for every $n\in \N$,
every complete surface immersion $f: F\to \E^n$, and every smooth, simple closed curve $\Gamma$ on $F$ of length $l$. If $\int_{\Gamma^{\beta l}}|A|^2\, d\mathrm{vol}_2 \leq c$,
then $|{\cal{K}}|(\Gamma)\geq\delta$.
\end{proposition}

\begin{remark}\label{rem42}
{\em More precisely, we will prove the following inequality for $\Gamma\subseteq F$ as above and for all $t>0$:
\begin{equation} \label{eq41}
|{\cal{K}}|(\Gamma) \geq \frac{3}{5}\left(\frac{4}{3}\pi - \int_{\Gamma^t} |A|^2\,d\mathrm{vol}_2\right) - \frac{8l}{5t}.
\end{equation}
From (\ref{eq41}) one easily concludes that Proposition \ref{prop41} holds, e.\,g. for $\delta(c) = \frac{1}{5}(\frac{4}{3}\pi-c)$ and
$\beta(c) = \frac{4}{\frac{4}{3}\pi-c}$. From the proof of (\ref{eq41}) it is clear that (\ref{eq41}) is not sharp.}
\end{remark}

The proof of Proposition \ref{prop41} is based on Fenchel's inequality for the total curvature of closed curves in $\E^n$, the Gau{\ss}-Bonnet formula, and the Bol-Fiala technique that provides area bounds for $\Gamma^t$ depending on integral bounds on the Gaussian curvature $K$. In the smooth case the ideas by G. Bol \cite{Bol41} and Fiala \cite{Fia41} were developed and extended by
P.\,Hartman \cite{Har64}, see also chapter 4 of the book \cite{SST03} and chapter 2 of the book \cite{BuZa88}. We use Hartman's results to treat the difficulties arising from the non-differentiability of the distance function $d^F(\cdot, \Gamma)$ from $\Gamma$.\\

Since the proof of inequality (\ref{eq41}) is a combination of several estimates we first give a rough outline. The starting point is Fenchel's inequality that implies
$$
\int_0^t |{\cal{K}}|(\partial\Gamma^\tau)\,d\tau + \int_{\Gamma^t}|A|\,d\mathrm{vol}_2 \geq 2\pi t
$$
provided $\partial\Gamma^t\neq \emptyset$. Now one would like to compare $|{\cal{K}}|(\partial\Gamma^\tau)$ to $|{\cal{K}}|(\Gamma)$, using the Gauss-Bonnet formula and the fact that
$\int_{\Gamma^\tau}|K| d\mathrm{vol}_2 \leq \frac{1}{2} \int_{\Gamma^\tau} |A|^2\,d\mathrm{vol}_2$ can be assumed to be small. Here one encounters two problems. The first problem is that the Euler characteristic of $\Gamma^\tau$ might be negative. The second problem is that the boundary terms in the Gau{\ss}-Bonnet formula involve the geodesic curvature, and not its absolute value. To overcome these problems we treat separately the signed total curvature ${\cal{K}}(\partial\Gamma^\tau)$ and its positive part ${\cal{K}}^{+}(\partial\Gamma^\tau)$. Since, for most
$\tau\in (0, t)$, $\frac{d}{d\tau}({\mathrm{length}}(\partial\Gamma^\tau))={\cal{K}}(\partial\Gamma^\tau)$, Hartman's results imply
$$
\int_0^t {\cal{K}}(\partial\Gamma^\tau)\,d\tau > -2l,
$$
see Lemma \ref{lemma44}. On the other hand, ${\cal{K}}^{+}(\partial\Gamma^\tau)$ can be bounded above by $|{\cal{K}}|(\Gamma)$ and $\frac{1}{2}\int_{\Gamma^\tau} |A|^2 d\mathrm{vol}_2$,
see Lemma \ref{lemma43}. This is a consequence of the Gau{\ss}-Bonnet formula and the fact that $\partial\Gamma^\tau$ has constant distance from $\Gamma$. Since
$|{\cal{K}}|(\partial\Gamma^\tau) = 2{\cal{K}}^{+}(\partial\Gamma^\tau) - {\cal{K}}(\partial\Gamma^\tau)$, the preceding inequalities combine to the inequality
$$
2|{\cal{K}}|(\Gamma)t \geq (2\pi - \int_{\Gamma^t}|A|^2\,d\mathrm{vol}_2)t - \int_{\Gamma^t}|A|\,d\mathrm{vol}_2 - 2l.
$$
Finally, we use $\int_{\Gamma^t}|A|\,d\mathrm{vol}_2 \leq (\int_{\Gamma^t}|A|^2\,d\mathrm{vol}_2)^{\frac{1}{2}} \mathrm{vol}_2(\Gamma^t)^{\frac{1}{2}}$, and the Bol-Fiala type estimate
$$
\mathrm{vol}_2(\Gamma^t)\leq 2lt + \left(\frac{1}{2}|{\cal{K}}|(\Gamma) + \frac{1}{4} \int_{\Gamma^t}|A|^2\,d\mathrm{vol}_2\right)t^2,
$$
see (\ref{eq49}), to obtain (\ref{eq41}).\\

Now we give the details of the proof of (\ref{eq41}). Let $\bar{R}=\sup_{x\in F} d^F(x, \Gamma)\in (0, \infty]$. Then $\partial\Gamma^t \neq \emptyset$ for $t\in (0, \bar{R})$. P.\,Hartman \cite{Har64} introduced the set ${\cal{NE}}={\cal{NE}}(\Gamma)\subseteq (0, \bar{R})$ of non-exeptional values of the distance function from $\Gamma$, and proved that ${\cal{NE}}$ is open
and of full measure in $(0, \bar{R})$. Note that if $t\in{\cal{NE}}$ then $\partial\Gamma^t$ is free of focal points of $\Gamma$, and for every $q\in\partial\Gamma^t$ there exist at most two geodesics of length $t$ joining $q$ to $\Gamma$. Moreover, if there are two such geodesics they intersect at $q$ at an angle smaller than $\pi$. This implies that $\partial\Gamma^t$ is a piecewise smooth submanifold of $F$, and that the set $Q^t\subseteq \partial\Gamma^t$ of points in the neighborhood of which $\partial\Gamma^t$ is not smooth, is finite. Moreover, $\Gamma^t$ has a concave angle at each $q\in Q^t$. For $t\in{\cal{NE}}$ we let $\kappa_g^t:\partial\Gamma^t\backslash Q^t \to \R$ denote the geodesic cuvature of $\partial\Gamma^t\backslash Q^t$
with respect to the normal pointing out of $\Gamma^t$, and $(\kappa_g^t)^{+}$ its positive part.

\begin{lemma}\label{lemma43}
If $t\in{\cal{NE}}$ then $\D\int_{\partial\Gamma^t\backslash Q^t}(\kappa_g^t)^{+}\, ds \leq |{\cal{K}}|(\Gamma)+ \int_{\Gamma^t}K^{-}\,d\mathrm{vol}_2$.
\end{lemma}

\proof Let $J \subseteq \partial\Gamma^t\backslash Q^t$ be a compact interval on which $\kappa_g^t$ is positive. The shortest connections to $\Gamma$ from the end-points of $J$,
together with $J$ and the nearest point projection ${\mathrm{pr}}(J)\subseteq\Gamma$ of $J$ to $\Gamma$ bound a rectangle $R_J$ in $\Gamma^t$. We choose an arclength-parametrization $\gamma$ of $\Gamma$, and we let $n_J$ denote the unit normal along ${\mathrm{pr}}(J)$ pointing into $R_J$. Then the Gau{\ss}-Bonnet formula yields
\begin{equation} \label{eq42}
\int_J(\kappa_g^t)^{+}(s)\,ds= \int_{{\mathrm{pr}}(J)}\langle \nabla_{\dot{\gamma}} \dot{\gamma}, n_J\rangle\,ds - \int_{R_J} K\,d\mathrm{vol}_2.
\end{equation}
If $J, J'$ are two such intervals and $J\cap J' = \emptyset$, then $\mathrm{vol}_2(R_J\cap R_{J'})=0$ and $n_J|{\mathrm{pr}}(J)\cap{\mathrm{pr}}(J')=
-n_{J'}|{\mathrm{pr}}(J)\cap{\mathrm{pr}}(J')$. Hence we have
\begin{equation} \label{eq43}
 \int_{{\mathrm{pr}}(J)\cap{\mathrm{pr}}(J')}\langle \nabla_{\dot{\gamma}} \dot{\gamma}, n_J\rangle\,ds =
-\int_{{\mathrm{pr}}(J)\cap{\mathrm{pr}}(J')}\langle \nabla_{\dot{\gamma}} \dot{\gamma}, n_{J'}\rangle\,ds.
\end{equation}
So, if we let $B\subseteq \Gamma$ denote the set of points $x\in \Gamma$ such that there exists precisely one $q\in\partial\Gamma^t\backslash Q^t$ with $\kappa_g^t(x)>0$ and
${\mathrm{pr}}(q)=x$, then (\ref{eq42}) and (\ref{eq43}) imply
$$
\int_{\partial\Gamma^t\backslash Q^t}(\kappa_g^t)^{+}\, ds \leq \int_B |\nabla_{\dot{\gamma}} \dot{\gamma}|\,ds + \int_{\Gamma^t} K^{-}\,d\mathrm{vol}_2
\leq |{\cal{K}}|(\Gamma)+\int_{\Gamma^t} K^{-}\,d\mathrm{vol}_2.
$$
\kasten

%% S. 4.5
If $t\in{\cal{NE}}(\Gamma) = {\cal{NE}}$ and $q\in Q^t \subseteq\partial\Gamma^t$, we let $\theta_q\in (0, \pi)$ denote the angle at $q$ between the two shortest connections from $q$ to $\Gamma$.
Then $\pi + \theta_q\in(\pi, 2\pi)$ is the inner angle of $\Gamma^t$ at $q$. The total curvature ${\cal{K}}(t)$ of $\partial\Gamma^t$ as the boundary of $\Gamma^t$ is defined by
$$
{\cal{K}}(t) = \int_{\partial\Gamma^t\backslash Q^t}\kappa_g^t\, ds - \sum_{q\in Q^t} \theta_q.
$$
\begin{lemma}\label{lemma44}
If $t\in{\cal{NE}}\subseteq (0, \bar{R})$, then $\int_0^t {\cal{K}}(\tau)\,d\tau > - 2l$, where $l={\mathrm{length}}(\Gamma)$.
\end{lemma}
Note. Since ${\cal{NE}}$ has full measure in $(0, \bar{R})$, ${\cal{K}}({\tau})$ is defined for almost all $\tau\in (0, \bar{R})$.\\

\proof
For $\tau\in{\cal{NE}}$  we let $l(\tau)$ denote the length of $\partial\Gamma^\tau$. Note that $\lim_{\tau\da 0} l(\tau) = 2l$.
A standard calculation shows that $l|{\cal{NE}}$ is smooth and that
\begin{equation} \label{eq44}
l'(\tau) = \int_{\partial\Gamma^\tau\backslash Q^\tau} \kappa_g^{\tau}\, ds - \sum_{q\in Q^\tau} 2 \tan \frac{\theta_q}{2},
\end{equation}
see e.\,g. \cite{SST03}, Theorem 4.1. Applying Hartman's Corollary 6.1 in \cite{Har64} we obtain
%it is a nontrivial consequence of Hartman's work that an 'integrated version' 
%of (\ref{eq44}) holds in form of the inequality
\begin{equation} \label{eq45}
l(t) \leq 2l +\int_0^t l'(\tau)\,d\tau \quad \mbox{ for all }t \in {\cal{NE}}.
\end{equation}
%for $t\in{\cal{NE}}$, see \cite{Har64}, Corollary 6.1. 
%Recall that $l'(\tau)$ is defined on the full measure subset ${\cal{NE}}$ 
%of $(0,\bar{R})$, and that $l'(\tau)$ satisfies (\ref{eq44}) for all $\tau \in {\cal{NE}}$.
Since $2\tan\frac{x}{2}\geq x$ for $x\geq 0$, we conclude from (\ref{eq44}) that 
$l'(\tau)\leq {\cal{K}}({\tau})$ for $\tau\in{\cal{NE}}$. 
Hence (\ref{eq45}) implies that
$$
l(t)\leq 2l+\int_0^t {\cal{K}}({\tau})\,d\tau \quad \mbox{ for }t \in {\cal{NE}}.
%\vspace{-0.4cm}
$$
This proves our claim since $l(t)> 0$ for $t \in {\cal{NE}}$.
\kasten

%% S. 4.6
\begin{lemma}\label{lemma45}
If $t\in (0, \infty)$ and $\int_{\Gamma^t} |A|^2\,d\mathrm{vol}_2 < 8\pi$, then $t< \bar{R}$ and $\partial\Gamma^t\neq \emptyset$.
\end{lemma}
\proof If $t\geq\bar{R}$, then $\Gamma^t = \Gamma^{\bar{R}}= F$, since $\bar{R}=\sup_{x\in F} d^F(x, \Gamma)$. In particular, $F$ is compact.
Then the well-known results by Chern-Lashoff \cite{CheLa57} imply that
$$
\int_F |A|^2\,d\mathrm{vol}_2 \geq 8\pi.\vspace{-0.4cm}
$$
\kasten

\textsc{Proof of inequality (\ref{eq41}):} Since (\ref{eq41}) is trivially true if $\int_{\Gamma^t} |A|^2\,d\mathrm{vol}_2\geq 8\pi$, we may assume that
$\int_{\Gamma^t} |A|^2\,d\mathrm{vol}_2 < 8\pi$, so that $t<\bar{R}$ by Lemma \ref{lemma45}. For $\tau\in{\cal{NE}}\cap (0, t)$ we use Lemma \ref{lemma43} to estimate
\begin{equation} \label{eq46}
\int_{\partial\Gamma^\tau\backslash Q^\tau} |\kappa_g^{\tau}|\,ds = \int_{\partial\Gamma^\tau\backslash Q^\tau} (2(\kappa_g^{\tau})^{+} - \kappa_g^{\tau})\,ds
\leq 2|{\cal{K}}|(\Gamma)+ 2\int_{\Gamma^\tau}K^{-}\,d\mathrm{vol}_2 - \int_{\partial\Gamma^\tau\backslash Q^\tau}\kappa_g^{\tau}\,ds.
\end{equation}
For $q\in \partial\Gamma^\tau\backslash Q^\tau$ we let $\kappa^{\tau}(q) \geq 0$ denote the curvature at $q$ of the curve $f|\partial\Gamma^\tau\backslash Q^\tau$ in $\E^n$.
Note that $\kappa^\tau(q) \leq |\kappa_g^{\tau}(q)|+ |A(q)|$. Hence Fenchel's inequality yields
\begin{equation} \label{eq47}
2\pi n(\tau) \leq \int_{\partial\Gamma^\tau\backslash Q^\tau} (|\kappa_g^{\tau}|+|A|)\,ds + \sum_{q\in Q^{\tau}} \theta_q,
\end{equation}
where $n(\tau)\geq 1$ denotes the number of connected components of $\partial\Gamma^\tau$.
Combining inequalities (\ref{eq46}), (\ref{eq47}) and $n(\tau)\geq 1$ we obtain
$$
2\pi \leq -{\cal{K}}(\tau) +2|{\cal{K}}|(\Gamma)+ 2\int_{\Gamma^\tau}K^{-}\,d\mathrm{vol}_2 + \int_{\partial\Gamma^\tau}|A|\,ds.
$$
Integrating this inequality over ${\cal{NE}}\cap (0, t)$, and using Lemma \ref{lemma44}, and $K^{-}\leq \frac{1}{2}|A|^2$, and the Cauchy-Schwarz inequality on
$\int_{\Gamma^t}|A|\,d\mathrm{vol}_2$, we obtain
\begin{equation} \label{eq48}
2\pi t \leq 2l + \big(2|{\cal{K}}|(\Gamma)+ \int_{\Gamma^t}|A|^2\,d\mathrm{vol}_2\big)t + (\int_{\Gamma^t}|A|^2\,d\mathrm{vol}_2)^{\frac{1}{2}} \cdot
\mathrm{vol}_2(\Gamma^t)^{\frac{1}{2}}
\end{equation}
To estimate $\mathrm{vol}_2(\Gamma^t)^{\frac{1}{2}}$ we recall that, by (\ref{eq44}) and Lemma \ref{lemma43},
$$
l'(\tau) \leq |{\cal{K}}|(\Gamma)+ \int_{\Gamma^\tau}K^{-}\,d\mathrm{vol}_2 \leq |{\cal{K}}|(\Gamma)+ \frac{1}{2}\int_{\Gamma^\tau} |A|^2\,d\mathrm{vol}_2.
$$
Combined with (\ref{eq45}) this implies
\begin{equation} \label{eq49}
\mathrm{vol}_2(\Gamma^{t}) = \int_0^t l(\tau)\,d\tau \leq 2lt + \frac{1}{2}\big(|{\cal{K}}|(\Gamma)+ \frac{1}{2}\int_{\Gamma^t}|A|^2\,d\mathrm{vol}_2\big)t^2
\end{equation}
and, hence,
$$
\mathrm{vol}_2(\Gamma^{t})^{\frac{1}{2}} \leq\left(\frac{1}{2}|{\cal{K}}|(\Gamma)+ \frac{1}{4}\int_{\Gamma^t}|A|^2\,d\mathrm{vol}_2\right)^{\frac{1}{2}}t +
l\left(\frac{1}{2}|{\cal{K}}|(\Gamma)+ \frac{1}{4}\int_{\Gamma^t}|A|^2\,d\mathrm{vol}_2\right)^{-\frac{1}{2}}.
$$
Inserting this into (\ref{eq48}) and abbreviating $\frac{1}{2}|{\cal{K}}|(\Gamma)=k$, $\frac{1}{4}\int_{\Gamma^t}|A|^2\,d\mathrm{vol}_2=a$, we obtain
$$
2\pi t \leq 2\big(1+\big(\frac{a}{k+a}\big)^{\frac{1}{2}}\big)l + (4(k+a)+2(a^2+ka)^{\frac{1}{2}})t\leq 4l + (5k+6a)t.
$$
This last inequality is equivalent to (\ref{eq41}). \kasten

\setcounter{equation}{0}
%%%%%%%%%%%%%%%%%%%%%%%%%%%%%%%%%%%%%%%%%%%%%%%%%%%%%%%%%%%%%%%
\section{Contractibility of intrinsic metric balls on surfaces in euclidean spaces} \label{sec5}

Using Proposition \ref{prop41} we will prove the following result that is the key to obtain a lower bound on the area of balls on surfaces in a compact Riemannian manifold, see Section \ref{sec6}.

\begin{proposition} \label{prop51}
Given $c\in \left(0, \frac{4}{3}\pi\right)$ there exists $\alpha= \alpha(c)\in(0,1)$ such that the following holds for every $n\in \N$, every complete surface immersion $f: F\to \E^n$, every $p\in F$ and every $R>0$. If $\int_{B(p, R)} |A|^2\,d\mathrm{vol}_2\leq c$, then the intrinsic metric ball $B(p, \alpha R)\subseteq F$ is simply connected.
\end{proposition}

\begin{remark}\label{rem52} {\em We do not know the supremum $c_0$ of all constants $c$ for which the statement in Proposition \ref{prop51} holds. By Proposition \ref{prop51} we have
$c_0\geq \frac{4}{3}\pi$, while the example of the catenoid shows that $c_0 \leq 8\pi$.}
\end{remark}

\begin{remark}\label{rem53} {\em Our proof provides an explicit value for $\alpha(c)$, namely}
$$
\alpha(c)=\left(\frac{2}{1-\cos\delta(c)}+ 2\beta(c)+1\right)^{-1},
$$
{\em where} $\delta(c)=\frac{1}{5}\left(\frac{4}{3}\pi-c\right)$ and $\beta(c)=\frac{4}{\frac{4}{3}\pi-c}$.
\end{remark}

The proof of Proposition \ref{prop51} depends crucially on Proposition \ref{prop41}. We will assume that $B(p, r)$ is not simply connected. Then we will prove in a series of lemmas that in some larger ball $B(p, R)$ we can find a simple geodesic loop $\Gamma$ that has an angle smaller than $\delta(c)$ at its base point, i.\,e. 
$|{\cal{K}}|(\Gamma)<\delta(c)$. This will contradict
Proposition \ref{prop41}. Here, the main idea is that the angle of a geodesic loop determines the speed at which it can be shortened by letting its base point vary. This is a purely intrinsic argument that does not depend on the assumption that $\dim F=2$. In the final step, however, we use the Gau{\ss}-Bonnet formula to see that 
$B(p, r)$ is not contractible in $B(p, R)$. Here and in the sequel a geodesic loop in a 
Riemannian manifold $M$ means a geodesic $\gamma:[0,L] \to M$ such that $\gamma(0) = \gamma(L)$,
and $\gamma(0)$ is called the base point of the loop. The loop will be called simple 
if $\gamma|_{[0,L)}$ is injective. 

\begin{lemma}\label{lemma54}
Let $\gamma: [0, l]\to M$ be a geodesic loop in a Riemannian manifold $M$ that is shortest among all non-contractible loops based at $\gamma(0)$. 
Then $\gamma$ is simple.
\end{lemma}

\proof Otherwise there exist $0\leq s_1 < s_2 < l$ such that $\gamma(s_1)= \gamma(s_2)$. Reversing the orientation of $\gamma$, if necessary, we may assume that $s_1\leq l-s_2$.
Note that the loop $\gamma|[s_1, s_2]$ is non-contractible since otherwise the loop $(\gamma|[0, s_1])\ast (\gamma|[s_2, l])$ were non-contractible, based at $\gamma(0)$, and shorter than $\gamma$. Hence the loop $\tilde{\gamma}= (\gamma|[0, s_2])\ast (\gamma|[0, s_1])^{-1}$ is non-contractible, based at $\gamma(0)$, and has length $s_2+s_1\leq l$. In particular, $\tilde{\gamma}$ is geodesic, in contradiction to the fact that $\dot{\gamma}(s_2)\neq - \dot{\gamma}(s_1)$.
\kasten\vspace{0.2cm}

\begin{lemma}\label{lemma55}
Let $B(p,r)$ be a metric ball in a complete Riemannian manifold, and suppose $B(p, r)$ is not simply connected. Then, among the arclength-parametrized loops in $B(p, r)$ with base point $p$
that are non-contractible in $B(p, r)$, there exists a shortest one, and every such shortest loop is a simple geodesic loop of length smaller than $2r$.
\end{lemma}

\proof The existence of a geodesic loop $\gamma:[0, l]\to B(p,r)$ that is shortest among all non-contractible loops in $B(p,r)$ based at $p$, and of length $l<2r$, is proved in \cite{Shio99}, Lemma 3.1. Lemma \ref{lemma54}, applied to $M=B(p,r)$, shows that every such loop is simple.
\kasten\vspace{0.2cm}

\begin{lemma}\label{lemma56} Suppose $\delta\in(0, \pi)$, $r>0$, and $R\geq \left(\frac{2}{1-\cos\delta}+1\right)r$. Let $B(p,R)$ be a metric ball of radius $R$ in a complete
Riemannian manifold, and assume that there exists a loop $\gamma_0$ based at $p$ 
of length smaller than $2r$ that is not contractible in $B(p,R)$. Then there exists a simple geodesic loop
$\gamma: [0, l]\to B(p,R)$ of length $l<2r$ such that the angle $\sphericalangle(\dot{\gamma}(0), \dot{\gamma}(l))$ that $\gamma$ makes at $\gamma(0) = \gamma(l)$ is smaller than $\delta$.
\end{lemma}

\proof For $\rho\in [0, R-r)$ we let $l(\rho)$ denote the infimum of the lengths of loops in $B(p,R)$ with base-points in $\overline{B(p, \rho)}$ that are not contractible in $B(p,R)$. Note that $l(\rho)$ is non-increasing and 2-Lipschitz. Since $l(0)<2r$ by assumption, the infimum $l(\rho)$ is achieved by a geodesic loop
$\gamma_\rho: [0, l(\rho)]\to B(p,\rho+r)\subseteq B(p,R)$. Lemma \ref{lemma54} implies that $\gamma_\rho$ is simple. Whenever $l'(\rho)$ exists, the first variation formula implies that $l'(\rho)\leq -1+\cos\delta_\rho$, where $\delta_\rho=\sphericalangle(\dot{\gamma}_\rho(0), \dot{\gamma}_\rho(l(\rho)))$. If, contrary to our claim, we had $\delta_\rho \geq \delta$ for all $\rho\in [0, R-r]$, we would have $0\leq l(R-r)<2r-(1-\cos\delta)(R-r)$, in contradiction to our assumption
$R\geq\left(\frac{2}{1-\cos\delta}+1\right)r$.
\kasten\vspace{0.2cm}

\textsc{Proof of Proposition \ref{prop51}:} Contrary to the claims of Proposition \ref{prop51} and Remark \ref{rem53} we assume that there exist $c\in(0, \frac{4}{3}\pi)$, a complete surface immersion $f:F\to \E^n$, $p\in F$ and $R>0$, such that $\int_{B(p, R)}|A|^2\,d\mathrm{vol}_2\leq c$, while $B(p, \alpha R)$ is not simply connected where
$\alpha= \left(\frac{2}{1-\cos\delta(c)}+ 2\beta(c)+1\right)^{-1}$. We set $r=\alpha R$ and apply Lemma \ref{lemma55} to obtain a simple geodesic loop
$\gamma_0: [0, l_0]\to B(p,r)$ of length $l_0<2r$ that is based at $p$ and not contractible in $B(p,r)$. Next we show that neither $\gamma_0$ is contractible in
$B(p,R)$. Otherwise $\gamma_0$ would bound a topological disk $D\subseteq B(p,R)$. Then, by the Gau{\ss}-Bonnet formula, we would have
$\int_D K\,d\mathrm{vol}_2= 2\pi-\sphericalangle(\dot{\gamma}_0(0), \dot{\gamma}_0(l_0))>\pi$, in contradiction to
$\int_D K\,d\mathrm{vol}_2\leq \frac{1}{2}\int_{B(p, R)} |A|^2 \,d\mathrm{vol}_2\leq \frac{1}{2}c < \pi$. Now we can apply Lemma \ref{lemma56} to $\delta(c)$, $r$ and
$\tilde{R}=\left(\frac{2}{1-\cos\delta(c)}+ 1\right)r<R$ to obtain a simple geodesic loop $\gamma:[0, \tilde{l}]\to B(p, \tilde{R})$ of length $\tilde{l}<2r$ such that
$\sphericalangle(\dot{\gamma}(0), \dot{\gamma}(\tilde{l}))<\delta(c)$. Smoothing the possible corner of $\gamma$ at $\gamma(0)=\gamma(\tilde{l})$, we obtain a smooth, simple closed curve $\Gamma\subseteq B(p, \tilde{R})$ of total absolute geodesic curvature $|{\cal{K}}|(\Gamma)<\delta(c)$ and of length $\tilde{l}<2r$. Note that, by the choice of $r=\alpha R$ and $\tilde{R}$, the tubular neighborhood $\Gamma^{\beta(c)l}$ of $\Gamma$ is contained in $B(p, R)$ so that the assumption
$\int_{\Gamma^{\beta(c) l}}|A|^2\, d\mathrm{vol}_2 \leq c$ in Proposition \ref{prop41} is satisfied. Hence the inequality $|{\cal{K}}|(\Gamma)<\delta(c)$ contradicts Proposition \ref{prop41}.
\kasten\vspace{0.2cm}

\setcounter{equation}{0}
%%%%%%%%%%%%%%%%%%%%%%%%%%%%%%%%%%%%%%%%%%%%%%%%%%%%%%%%%%%%%%%
\section{Area of intrinsic metric balls on surfaces in Riemannian manifolds} \label{sec6}

We consider a complete surface immersion $f:F\to M$ into a Riemannian manifold $M$. We assume that the sectional curvature of $M$ is bounded below by some negative number $k<0$.
We will derive upper and lower estimates for the area of metric balls on $F$, that depend on assumptions on the energy. The upper estimate is standard, and has its origin in an idea that goes back to work by G.\, Bol \cite{Bol41} and by F.\, Fiala \cite{Fia41}:

\begin{proposition} \label{prop61}
Let $F$ be a complete Riemannian surface, and let $K$ denote its Gaussian curvature. Then the following holds for every number $k<0$, for every $p\in F$ and every $r>0$:
$$
\mathrm{vol}_2\big(B(p, r)\big) \leq \frac{1}{-k}\big(2\pi + \int_{B(p, r)}(K-k)^{-}\,d \mathrm{vol}_2\big)\big(\cosh(\sqrt{-k}\,r)-1\big).
$$
\end{proposition}

\proof This follows easily from case 1) of Theorem 2.4.2 in \cite{BuZa88}. Note also that Corollary \ref{coroll115} will generalize Proposition \ref{prop61}.
\kasten\vspace{0.2cm}

\begin{corollary}\label{coroll62} Let $f:F\to M$ be a complete surface immersion into a Riemannian manifold $M$ with sectional curvature bounded below by $k<0$.
Then the following holds for every $p\in F$ and every $r>0$:
$$
\mathrm{vol}_2(B(p, r)) \leq \frac{1}{-k} (2\pi +E(f|B(p, r)))\big(\cosh(\sqrt{-k}\,r)-1\big).
$$
\end{corollary}
\proof
By the Gau{\ss} equations we have $K\geq k-\frac{1}{2}|A|^2$, cf. (1.2). Hence we have $(K-k)^{-}\leq \frac{1}{2}|A|^2$ and
$\int_{B(p, r)}(K-k)^{-}\,d \mathrm{vol}_2\leq E(f|B(p, r))$. So our claim is reduced to Proposition \ref{prop61}.
\vspace{0.2cm}\kasten

Our lower area estimate depends on the result of Section \ref{sec5}.

\begin{proposition} \label{prop63}
Let $M$ be a compact Riemannian manifold. Then there exist $r_0>0$ and $\beta\geq 1$ such that the following holds for every complete surface immersion $f:F\to M$ and every $p\in F$.
If $r\in(0, r_0)$ and $E(f|B(p, \beta r))\leq \frac{\pi}{4}$, then $B(p, r)$ is simply connected, and $\mathrm{vol}_2(B(p,r))\geq\frac{\pi}{2}r^2$.
\end{proposition}

\proof
Let $i:M\to\E^n$ be an isometric immersion of $M$ into some euclidean space, and let $C>0$ be an upper bound for the norm of the second fundamental form
$\tilde{A}$ of $i$. Then the second fundamental form $\bar{A}$ of $\bar{f}=i\circ f$ can be orthogonally decomposed as $\bar{A}=A + f^{\ast}\tilde{A}$,
where $A$ denotes the second fundamental form of $f$. This implies $|\bar{A}|^2\leq |A|^2+C^2$. Hence, using Corollary \ref{coroll62} we obtain
$$
E(\bar{f}|B(p,r))\leq E(f|B(p, r)) + \frac{C^2}{-2k}(2\pi+E(f|B(p, r)))(\cosh(\sqrt{-k}\,r)-1),
$$
where $k<0$ denotes a lower bound for the sectional curvature of $M$. According to the Gau{\ss} equations we can take $k=-C^2$, so that the preceding inequality reduces to
$$
E(\bar{f}|B(p,r))\leq \pi(\cosh(Cr)-1)+\frac{1}{2}E(f|B(p, r))(\cosh(Cr)+1).
$$
Hence we can find $r_1>0$ so that $r\in (0, r_1)$ and $E(f|B(p, r))\leq \frac{\pi}{4}$ imply
\begin{equation}\label{eq61}
E(\bar{f}|B(p, r))\leq \frac{\pi}{2},
\end{equation}
namely $r_1= \frac{1}{C}\arcosh(\frac{11}{9})$. Now we can apply Proposition \ref{prop51} to find $\alpha=\alpha(\pi)$ so that $B(p, \alpha r)$ is simply connected whenever $r\in(0, r_1)$ and $E(f|B(p, r))\leq \frac{\pi}{4}$. We set $\beta = \alpha^{-1}$ and $r_0=\alpha r_1$, and conclude that $B(p, r)$ is simply connected whenever $r\in(0, r_0)$ and $E(f|B(p, \beta r))\leq \frac{\pi}{4}$. This proves our first claim. For these $r$ we can apply Proposition 3.2(3) in \cite{Shio99}, and obtain
\begin{equation}\label{eq62}
\mathrm{vol}_2(B(p, r)) \geq (\pi-\int_{B(p, r)}K^{+}\,d \mathrm{vol}_2)r^2.
\end{equation}
If $r\in (0, r_0)$ and $E(f|B(p, \beta r))\leq \frac{\pi}{4}$, we use $K^{+}\leq \frac{1}{2}|\bar{A}|^2$ and (\ref{eq61}) to see that
$\int_{B(p, r)}K^{+}\,d \mathrm{vol}_2\leq E(\bar{f}|B(p,r))\leq \frac{\pi}{2}$. Hence (\ref{eq62}) implies
$\mathrm{vol}_2(B(p, r)) \geq \frac{\pi}{2}r^2$.
\kasten\vspace{0.2cm}

\setcounter{equation}{0}
%%%%%%%%%%%%%%%%%%%%%%%%%%%%%%%%%%%%%%%%%%%%%%%%%%%%%%%%%%%%%%%
\section{Gromov-Hausdorff convergence} \label{sec7}

We consider a sequence of complete surface immersions $f_i:F_i\to M$ into a compact Riemannian manifold $(M, g)$.
We assume that the $F_i$ and $M$ are connected, and let $d_i$ resp. $d^M$ denote the distance functions induced by $f_i^{\ast} g$ resp. by $g$.
Under appropriate assumptions on the energy of the $f_i$ we will study pointed Gromov-Hausdorff convergence, abbreviated as GH-convergence,
of the metric spaces $(F_i, d_i)$ to a limit metric space. We will rely on the concept of GH-convergence as explained by M.\,Gromov in \cite{Gro81}, Section 6.
Except for some technical details much of the material in this section is a special case of \cite{Shio99}, Section 3. Here, we present complete details since, on the one hand, the results from Section \ref{sec6} make the proofs in our special case considerably simpler than the ones given in \cite{Shio99}, while, on the other hand,
the results in \cite{Shio99} do not apply directly in our situation.\\

We start by recalling some notions from metric geometry. A metric space $X$ is {\em proper} if the compact subsets of $X$ are precisely the closed and bounded subsets of $X$. $X$ is a {\em length space} if any two points $x, y\in X$ can be joined by a curve in $X$ of length $d(x,y)$. A sequence of compact metric spaces $X_i$ is
{\em uniformly compact} if the diameters of the $X_i$ are uniformly bounded, and if, for every $\ve >0$, there exists $N\in \N$ such that each $X_i$ has an
$\ve$-dense subset with at most $N$ elements. We recall Gromov's compactness criterion \cite{Gro81}, p.\,64:\\

{\em Let $(X_i, x_i)_{i\in\N}$ be a sequence of proper, pointed metric spaces. If for each $r>0$ the sequence of balls $(\overline{B(x_i, r)})_{i\in\N}$ is uniformly compact, then a subsequence of the sequence $(X_i, x_i)_{i\in\N}$ GH-converges to a proper, pointed metric space.}\\

Finally, we note the following well-known fact, see e.\,g. \cite{BBI01}, Theorem 7.5.1:
\begin{remark} \label{rem71}
{\em If the sequence $(X_i, x_i)_{i\in\N}$ GH-converges to a complete, pointed metric space $Y$, and if all the $X_i$ are length spaces, then so is $Y$.}
\end{remark}

If $X$ is a metric space and $\ve >0$, we define $\alpha_X(\ve)\in\N\cup\{\infty\}$ by
\begin{equation}\label{eq71}
\alpha_X(\ve)= \mbox{minimal number of elements of an $\ve$-dense subset of } X.
\end{equation}
If $A$ is a subset of $X$, we consider $A$ with the metric induced from $X$. This defines $\alpha_A(\ve)$ for subsets $A$ of $X$.\\

The following lemma relies on the area estimates obtained in Section \ref{sec6}. In particular, the statement involves the constants $r_0>0$ and $\beta\geq 1$ from
Proposition \ref{prop63}. As in Section \ref{sec6}, we let $k<0$ denote a lower bound for the sectional curvature of $M$.

\begin{lemma}\label{lemma72}
Let $f:F\to M$ be a complete surface immersion, and suppose $p\in F$ and $\tilde{R}>0$ are such that $E(f|B(p, \tilde{R}))\leq \frac{\pi}{4}$. If $0<R<\tilde{R}$
and $0 <\ve < 2{\mathrm{min}}\{r_0, \frac{\tilde{R}-R}{\beta}\}$, then
$$
{\alpha\,}_{\overline{B(p, R)}}(\ve)\leq \frac{18}{-k}(\cosh(\sqrt{-k}\tilde{R})-1)\ve^{-2}.
$$
\end{lemma}

\proof From Corollary \ref{coroll62} we know that
\begin{equation}\label{eq72}
\mathrm{vol}_2(B(p, \tilde{R})) \leq \frac{9\pi}{-4k}\big(\cosh(\sqrt{-k}\,\tilde{R})-1\big).
\end{equation}
If $D\subseteq \overline{B(p, R)}$ is $\ve$-dense in $\overline{B(p, R)}$ and $\# D ={\alpha}_{\overline{B(p, R)}}(\ve)$, then
$B(x, \frac{\ve}{2})\cap B(y, \frac{\ve}{2})=\emptyset$ for $x\neq y$ in $D$. Since $\frac{\ve}{2}< \tilde{R}-R$ we have
$B(x, \frac{\ve}{2})\subseteq B(p, \tilde{R})$ for all $x\in D$. This implies
\begin{equation}\label{eq73}
\mathrm{vol}_2(B(p, \tilde{R})) \geq \sum_{x\in D} \mathrm{vol}_2(B(x, \frac{\ve}{2})).
\end{equation}
Finally, we have $ \frac{\ve}{2}< r_0$ and $E(f|B(x, \beta\frac{\ve}{2}))\leq\frac{\pi}{4}$ for $x\in D$, since
$B(x, \beta\frac{\ve}{2})\subseteq B (p, R+\beta\frac{\ve}{2})\subseteq B(p, \tilde{R})$.
So we can apply Proposition \ref{prop63} to obtain
\begin{equation}\label{eq74}
\mathrm{vol}_2(B(x,\frac{\ve}{2}))\geq \frac{\pi}{8} \ve^2\quad\mbox{if } x\in D.
\end{equation}
Now inequalities (\ref{eq72})--(\ref{eq74}) imply our claim.
\kasten\vspace{0.2cm}

Using Gromov's compactness criterion and Lemma \ref{lemma72} we obtain:

\begin{proposition} \label{prop73}
\begin{itemize}
\item[\em{(a)}] Suppose $0<R<\tilde{R}$ and $p_i\in F_i$ is a sequence such that $\liminf_{i\to\infty} E(f_i|B(p_i, \tilde{R}))<\frac{\pi}{4}$. Then a subsequence of the sequence
$(\overline{B(p_i, R)}, p_i)_{i\in \N}$ GH-converges to a connected, compact, pointed metric space $(Y_R, y_0)$ of finite 2-dimensional Hausdorff measure.
\item[\em{(b)}] Suppose $p_i\in F_i$ is a sequence such that $\liminf_{i\to\infty} E(f_i|B(p_i, R))<\frac{\pi}{4}$ holds for all $R>0$. Then a subsequence of the sequence
$(F_i, p_i)_{i\in \N}$ GH-converges to a proper, pointed length space $(Y, y_0)$ of locally finite 2-dimensional Hausdorff measure.
\end{itemize}
\end{proposition}

\proof Lemma \ref{lemma72} shows that the sequence $(\overline{B(p_i, R)})_{i\in \N}$ is uniformly compact, in case (a) for $R<\tilde{R}$, and in case (b) for all $R>0$.
So the statements concerning GH-convergence follow from Gromov's compactness criterion. The estimate ${\alpha\,}_{\overline{B(p_i, R)}}(\ve)\leq c(k, R)\ve^{-2}$ from Lemma \ref{lemma72}
implies a similar estimate for the limit spaces $Y_R$ resp. $Y$. This proves the 
statements about the 2-dimensional Hausdorff measure, see Section 2.3 and
the proof of Theorem 3.1 in \cite{Shio99}. Remark \ref{rem71} shows that $Y$ 
is a length space.
\kasten
%\vspace{0.2cm}
\begin{remark} {\em In \cite{Shio99} T. Shioya studies the topological structure of 
{GH}-limit spaces of sequences of compact Riemannian surfaces with bounds 
on the total absolute curvature and the diameter. These results do not
apply directly in the situation of Proposition \ref{prop73}, since the
assumptions are not satisfied. In Section 13 we will prove that, under 
the conditions (\ref{eq131})-(\ref{eq133}), the limit space $Y$ in 
Proposition \ref{prop73}(b) admits a locally isometric map onto a complete,
totally geodesic surface in $M$, see Proposition \ref{prop133}.}
\end{remark}

To define a notion of convergence for the sequence of maps $f_i$ we resort to {\em definite GH-convergence}, as defined in \cite{Gro81}, p.\,66. If $(X_i, x_i)_{i\in \N}$
GH-converges to $(Y, y_0)$ we can choose and will fix a sequence of metrics $\delta_i$ on the disjoint unions $X_i\cup Y$ with the following properties:
\begin{itemize}
\item[(a)] The inclusions $(X_i, d_i)\to (X_i\cup Y, \delta_i)$ and $(Y, d_Y)\to (X_i\cup Y, \delta_i)$ are isometric.
\item[(b)] For every $r>0$ and every $\ve>0$ the following holds for almost all $i\in \N$:\\
           $\delta_i(x_i, y_0)<\ve$, and $B(x_i, r)$ is contained in the $\ve$-neighborhood of $Y$, and $B(y_0, r)$ is contained in the $\ve$-neighborhood of $X_i$
           (both neighborhoods with respect to $\delta_i$).
\end{itemize}

\begin{definition} \label{def74} Suppose $(X_i, \bar{x}_i)_{i\in \N} \to (Y, y_0)$ with respect to definite GH-convergence. Then we define:
\begin{itemize}
\item[\em{(a)}] A sequence $x_i\in X_i$ converges to $y\in Y$, i.\,e. $\lim_{i\to\infty} x_i=y$, if $\lim_{i\to\infty}\delta_i(x_i, y)=0$.
\item[\em{(b)}] Let $f_i: X_i\to Z$ be a sequence of maps into a metric space $Z$. Then the $f_i$ converge to a map $f: Y\to Z$ if $\lim_{i\to\infty} f_i(x_i)=f(y)$,
whenever the sequence $x_i\in X_i$ converges to $y\in Y$.
\end{itemize}
\end{definition}

Slightly generalizing the usual proof of the Arzel\`a-Ascoli theorem, we obtain:
\begin{proposition} \label{prop75}
{\em{(a)}} Suppose the sequence $(\overline{B(p_i, R)}, p_i)$  converges to $(Y_R, y_0)$ with respect to definite GH-convergence. Then a subsequence of the sequence
$(f_i|\overline{B(p_i, R)})_{i\in\N}$ converges to a distance-nonincreasing map $f: Y_R\to M$.\\
{\em{(b)}} Suppose the sequence $(F_i, p_i)_{i\in\N}$ converges to $(Y, y_0)$ with respect to definite GH-convergence. Then a subsequence of the sequence
$(f_i)_{i\in\N}$ converges to a distance-nonincreasing map $f: Y\to M$.\\
\end{proposition}
Similarly, the Arzel\`a-Ascoli theorem implies:
\begin{proposition} \label{prop76}
Suppose the sequence $(\overline{B(p_i, R)}, p_i)_{i\in\N}$ converges to $(Y_R, y_0)$ and $\lim_{i\to\infty} (f_i|\overline{B(p_i, R)}) = f: Y_R\to M$.
Let $L>0$, and let $w_i\in SF_i$ be a sequence such that $c_{w_i}([-L, L])\subseteq \overline{B(p_i, R)}$ for all $i\in \N$, and
$\lim_{i\to\infty} \int_{-L}^{L} |A_i|^2\circ c_{w_i}(t)\,dt=0$, and $\lim_{i\to\infty} df_i(w_i)=w\in SM$. Then a subsequence of the
$(c_{w_i}|[-L, L])_{i\in\N}$ converges to a curve in $Y_R$, and every limit curve $c:[-L, L]\to Y_R$ satisfies $f\circ c = c_w|[-L, L]$
and $d_{Y_R}(c(s), c(t))=|t-s|$ for all $s, t\in [-L, L]$ with $|t-s|\leq {\mathrm{injrad}}(M)$. In particular, we have $c_w([-L,L])\subseteq f(Y_R)$.
\end{proposition}

\proof Since the GH-convergence of $(\overline{B(p_i, R)}, p_i)_{i\in\N}$ to $(Y_R, y_0)$ can be realized by Hausdorff convergence within some fixed compact metric space, cf. \cite{Gro81}, p.\,65, the usual Arzel\`a-Ascoli theorem implies the existence of a limit curve $c:[-L, L]\to Y_R$, say
$\lim_{i\to\infty} c_{w_i}|[-L, L]=c$. So, if $t_i\in [-L, L]$ and $\lim_{i\to\infty} t_i=t$, then $\lim_{i\to\infty} c_{w_i}(t_i)=c(t)$, while
$\lim_{i\to\infty} f_i\circ c_{w_i}(t_i)=c_w(t)$ by Proposition \ref{prop23}. From Definition \ref{def74}(b) we see that this implies $f_i\circ c(t)=c_w(t)$.
Since $c$ and $f$ are 1-Lipschitz we conclude that $d_{Y_R}(c(s), c(t))=|t-s|$ if $s, t\in [-L,L]$ and $|t-s|\leq {\mathrm{injrad}}(M)$.
\kasten\vspace{0.2cm}

\setcounter{equation}{0}
%%%%%%%%%%%%%%%%%%%%%%%%%%%%%%%%%%%%%%%%%%%%%%%%%%%%%%%%%%%%%%%
\section{Asymptotic density of good points} \label{sec8}

We continue to consider a sequence of complete surface immersions $f_i:F_i\to M$ into a compact Riemannian manifold $M$. We assume that there exist $\tilde{R}>0$ and a sequence $p_i\in F_i$ such that $\lim_{i\to\infty} E(f_i|B(p_i, \tilde{R}))=0$. Relying on the volume estimate from Proposition \ref{prop36} and the lower estimate for the area of balls from Proposition \ref{prop63}, we will prove that there exists a sequence $\ve_i\da 0$ such that the sets of ``$(|A_i|^2, \ve_i, R, r)$-good points'' $G_{\ve_i, R, r}(|A_i|^2)$, cf. Definition \ref{def35}, are asymptotically dense in $B(p_i, r)$, whenever $R>0$, $r>0$ and $R+2r<\tilde{R}$. In combination with the results from Section \ref{sec7} this will be used to prove the existence of a totally geodesic surface in $M$. As a prerequisite for the proof of a {\em complete} totally geodesic surface in $M$, we will then assume the existence of a sequence $\tilde{R}_i\to \infty$ such that $E(f_i|B(p_i, \tilde{R}_i))\to 0$ and prove the existence of sequences $\ve_i\da 0$, $R_i\to \infty$ such that the sets $\tilde{G}_{\ve_i, R_i}(|A_i|^2)$, cf. (\ref{eq37}), are asymptotically dense in
$B(p_i, R_i)$.
\begin{definition} \label{def81}
Let $(X_i, d_i)$ be a sequence of metric spaces, and, for each $i\in\N$, let $B_i$ and $G_i$ be subsets of $X_i$. Then the sequence $G_i$ is called asymptotically dense in the sequence
$B_i$ if
$$
\lim_{i\to\infty} (\sup\{d_i(x, G_i)|x\in B_i\})=0.
$$
\end{definition}

Put differently, this condition says that there is a sequence $\delta_i\da 0$ such that each $B_i$ is contained in the $\delta_i$-neighborhood of $G_i$. In the following we will depend on the constants $r_0>0$, $\beta\geq 1$ from Proposition \ref{prop63}.
\begin{lemma}\label{lemma82}
Assume there exist $\tilde{R}>0$ and a sequence $p_i\in F_i$ such that $\lim_{i\to\infty} E(f_i|B(p_i, \tilde{R})) = 0$. Then there exists a sequence $\ve_i\da 0$ such that the sequence $G_{\ve_i,R,r}(|A_i|^2)$ is asymptotically dense in the sequence $B(p_i, r)$, whenever $R>0$, $r>0$ and $R+2r<\tilde{R}$.
\end{lemma}

\proof We choose a sequence $\ve_i\da 0$ such that $\lim_{i\to\infty}(\ve_i^{-4} E(f_i|B(p_i, \tilde{R}))) = 0$. Since $R+2r<\tilde{R}$ we can use Proposition \ref{prop36}(b) to find $\tilde{r}>r$ such that
\begin{equation}\label{eq81}
\lim_{i\to\infty}(\mathrm{vol}_2(B(p_i, \tilde{r})\backslash G_{\ve_i,R,r}(|A_i|^2)))=0.
\end{equation}
If our claim were not true we could find $\ve>0$ and an infinite set $I\subseteq\N$ such that, for each $i\in I$, there exists $q_i\in B(p_i, r)$ with
$B(q_i, \ve)\cap G_{\ve_i,R,r}(|A_i|^2)=\emptyset$. We may assume that $\ve< \min\{r_0, \tilde{r}-r\}$ and $r+\beta\ve < \tilde{R}$.
Then $B(q_i, \beta\ve)\subseteq B(p_i, \tilde{R})$ and, hence, $E(f_i|B(q_i,\beta\ve))\leq \frac{\pi}{4}$ for almost all $i\in I$. For these infinitely many $i\in I$
Proposition \ref{prop63} implies $\volz(B(q_i,\ve))\geq \frac{\pi}{2} \ve^2$. On the other hand, since $\ve < \tilde{r}-r$, we have
$B(q_i,\ve)\subseteq B(p_i, \tilde{r})\backslash G_{\ve_i,R,r}(|A_i|^2)$, in contradiction to (\ref{eq81}).
\kasten

Similarly, we can prove:
\begin{lemma}\label{lemma83}
Assume there exist sequences $\tilde{R}_i\to \infty$ and $p_i\in F_i$ such that $\lim_{i\to\infty} E(f_i|B(p_i, \tilde{R}_i)) = 0$.
Then there exist sequences $\ve_i\da 0$, $R_i\to \infty$ such that the sequence $\tilde{G}_{\ve_i,R_i}(|A_i|^2)$ is asymptotically dense in the sequence $B(p_i, R_i)$.
\end{lemma}

\proof We can find sequences $\ve_i\da 0$ and $R_i\to \infty$ such that $5R_i<\tilde{R}_i$ for all $i\in\N$ and
$\lim_{i\to\infty} \left(\frac{R_i^3}{\ve_i^4} E(f_i|B(p_i, \tilde{R}_i))\right) = 0$. Then we can use Corollary \ref{coroll37} to conclude that
$$
\lim_{i\to\infty}(\mathrm{vol}_2(B(p_i, R_i+1)\backslash \tilde{G}_{\ve_i,R_i}(|A_i|^2)))=0.
$$
This implies the asymptotic density of $\tilde{G}_{\ve_i,R_i}(|A_i|^2)$ in $B(p_i, R_i)$ as in the proof of Lemma \ref{lemma82}.
\kasten

\setcounter{equation}{0}
%%%%%%%%%%%%%%%%%%%%%%%%%%%%%%%%%%%%%%%%%%%%%%%%%%%%%%%%%%%%%%%
\section{Existence of totally geodesic surfaces and proof of Theorem \ref{thrm12}} \label{sec9}

In this section we prove Theorem \ref{thrm12}. We assume that there exists a sequence of complete surface immersions $f_i:F_i\to M$ with points $p_i\in F_i$ and $\tilde{R}>0$ such that $\lim_{i\to\infty} E(f_i|B(p_i, \tilde{R})) = 0$, and prove that $M$ contains (a piece of) totally geodesic surface.\\

We start with a rough outline of the proof. According to Lemma \ref{lemma82} we may assume that $p_i\in G_{\ve_i,R,r}(|A_i|^2)$, where
$\ve_i\da 0$ and $R+2r< \tilde{R}$. Since $M$ is compact we may assume that the sequence $(df_i(T_{p_i}F_i))_{i\in\N}$ converges in the Grassmann bundle
$\pi_G:G_2M \to M$ of 2-dimensional linear subspaces in $TM$, say
$\lim_{i\to\infty} df_i(T_{p_i}F_i))=P\in (G_2 M)_p$, where $p=\lim_{i\to\infty}f_i (p_i)$.
We intend to prove that $N_P(R)=N=\{c_v(t)|v\in SP, |t|<R\}\subseteq M$ is a totally geodesic, embedded surface if $R$ is smaller than the injectivity radius of $M$ at $p$.
According to Propositions \ref{prop73} and \ref{prop75} we may assume that the sequence of compact metric spaces $(\overline{B(p_i, R+r)}, p_i)_{i\in \N}$ converges to a compact metric space
$(Y, y_0)$, and that the maps $f_i|B(p_i, R+r)$ converge to a 1-Lipschitz map $f:Y \to M$. Since $Y$ has finite 2-dimensional Hausdorff measure, the same is true for $f(Y)$.
Now, for fixed $w\in SP$, we define maps $\tilde{w}: SP\times(0, R)\to SM$, $\tilde{w}(v, t) = {\cal{P}}_t^{c_v}(w)$, and
$\psi^w=\psi: SP\times(0, R)\times(-r, r)\to M$, $\psi(v, t, s) = c_{\tilde{w}(v, t)}(s)$.
Then $\psi(SP\times(0, R)\times\{0\})= N\backslash\{p\}$. Using our assumption $p_i\in G_{\ve_i,R,r}(|A_i|^2)$ and Proposition \ref{prop76} we prove that
$\psi(SP\times(0, R)\times(-r,r))\subseteq f(Y)$. In particular, we see that ${\rm rk}(d\psi)= 2$ in a neighborhood of $SP\times(0, R)\times\{0\}$. Then standard calculus implies that
$\psi(v, t, s)=c_{\tilde{w}(v, t)}(s)\in \psi(SP\times(0, R)\times\{0\})\subseteq N$ for small enough $|s|$ and $\tilde{w}(v, t)\in S_{c_v(t)}N$. Since $\tilde{w}(v, t)$ is an arbitrary
vector in $S_{c_v(t)}N$, this proves that $N$ is totally geodesic.\\

Here are the precise formulation and the proof.

\begin{proposition} \label{prop91}
Let $f_i: F_i\to M$ be a sequence of complete surface immersions into a compact Riemannian manifold. Suppose $\ve_i\da 0$, $R>0$, $r>0$, and
$p_i\in G_{\ve_i, R,r}(|A_i|^2)$ is a sequence such that $\lim_{i\to\infty} E(f_i|B(p_i, R+2r)) = 0$ and $\lim_{i\to\infty} df_i(T_{p_i}F_i)=P\in (G_2 M)_p$.
If $R$ is smaller than the injectivity radius of $M$ at $p$, then $N_P(R) = N = \{c_v(t)|v\in SP, |t|<R\}\subseteq M$ is a totally geodesic, embedded surface.
\end{proposition}

\begin{remark} \label{rem92}
{\em The condition ``$R$ smaller than the injectivity radius of $M$ at $p$'' is imposed in order to have a convenient description of the totally geodesic surface $N$. The methods used in Section \ref{sec10} indicate that without this condition one should be able to find a complete surface immersion $f:F\to M$ and $q\in F$ such that $f|B(q, R)$ is totally geodesic.}
\end{remark}

\proof
According to Propositions \ref{prop73} and \ref{prop75} we can find a subsequence such that the sequence $(\overline{B(p_i, R+r)}, p_i)_{i\in \N}$ converges to a compact, pointed metric space
$(Y, y_0)$ with respect to definite GH-convergence and such that the sequence $(f_i|\overline{B(p_i, R+r)})_{i\in \N}$ %$(f_i|B(p_i, R+r))_{i\in \N}$
converges to a 1-Lipschitz map $f:Y \to M$. We fix an arbitrary $w\in SP$ and define maps $\tilde{w}:SP\times(0, R)\to TM$,
$\tilde{w}(v, t) = {\cal{P}}_t^{c_v}(w)$, and $\psi: SP\times(0, R)\times(-r, r)\to M$, $\psi(v, t, s)=c_{\tilde{w}(v, t)}(s)$, as above.
Since $R$ is smaller than the injectivity radius of $M$ at $p$, we know that $N$ is an embedded submanifold and $\psi|SP\times(0, R)\times\{0\}$ is a diffeomorphism onto $N\backslash \{p\}$.
We want to prove that $\tilde{w}(v, t)\in S_{c_v(t)}N$ and that $\psi(SP\times(0, R)\times(-r, r))\subseteq f(Y)$. Since $p_i\in G_{\ve_i,R,r}(|A_i|^2)$ we know that the sets
$(S_{p_i} F_i \times S_{p_i}F_i \times [-R, R])\backslash {\cal{V}}^2_{\ve_i, R, r}(|A_i|^2)$ are asymptotically dense in $S_{p_i} F_i \times S_{p_i}F_i \times [-R, R]$, cf. (\ref{eq34})
and Definition \ref{def35}. Hence, given $(v, w, t)\in SP\times SP\times (0, R)$, we can find a sequence $(v_i, w_i, t_i)\in S_{p_i} F_i \times S_{p_i}F_i \times (0, R)$ such that
$$
\lim_{i\to\infty} (df_i(v_i), df_i(w_i), t_i)= (v, w, t),
$$
and such that $\int_{-R}^R|A_i|^2\circ c_{v_i}(t)\,dt < \ve_i$ and $\int_{-r}^r |A_i|^2\circ c_{\tilde{w}_i}(s)\,ds < \ve_i$, where $\tilde{w}_i = {\cal{P}}_{t_i}^{c_{v_i}}(w_i)$, cf. (\ref{eq31}).
Now Proposition \ref{prop23} implies that
$$
\lim_{i\to\infty} df_i(\tilde{w}_i)= {\cal{P}}_t^{c_v}(w) = \tilde{w}(v, t),
$$
and Proposition \ref{prop76} implies that
$$
\psi(v, t, s) = c_{\tilde{w}(v, t)}(s)\in f(Y)
$$
if $s\in (-r, r)$. Since $f(Y)$ has finite 2-dimensional Hausdorff measure we conclude that 
${\mathrm{rk}}(d\psi)\leq 2$. Since ${\mathrm{rk}}(d\psi(v,t,0))\geq 2$ for
$(v, t)\in SP\times (0, R)$, we can find a neighborhood of $SP\times (0, R)\times\{0\}$ in $SP\times \R^2$ on which $d\psi$ has constant rank 2.
Using the local normal form of maps of constant rank we find, for every $(v, t)\in SP \times (0, R)$, some $\delta >0$ such that
$\psi(v,t,s)\in \psi(SP\times (0, R)\times\{0\})=N\backslash \{p\}$ for $|s|<\delta$. Since $N$ is an embedded submanifold this implies that
$\tilde{w}(v, t) = \frac{\partial \psi}{\partial s}(v, t, 0) \in S_{c_v(t)}N$. Hence ${\cal{P}}_t^{c_v}(SP)= S_{c_v(t)}N$, and the preceding argument shows that every geodesic
with initial vector in $S_{c_v(t)}N$ lies locally in $N$. This proves that $N$ is totally geodesic.
\kasten\vspace{0.2cm}

\begin{remark} \label{rem93}
{\em The preceding proof shows additionally that under the assumptions of Proposition \ref{prop91} every $q\in N$ is the limit of a sequence $(f_i(q_i))_{i\in \N}$
with $q_i\in B(p_i, R)\subseteq F_i$.\\

Theorem \ref{thrm12} is a consequence of the following slightly stronger result.}
\end{remark}

\begin{theorem}\label{thrm94}
Let $f_i: F_i \to M$ be a sequence of complete surface immersions into a compact Riemannian manifold. Assume the existence of a sequence $p_i\in F_i$ and of $\tilde{R}>0$ such that
$\lim_{i\to\infty} E(f_i|B(p_i, \tilde{R})) = 0$. Then there exists a 2-plane $P\in G_2 M$ such that $N_P(R)=\{c_v(t)|v\in SP, |t|<R\}$ is a totally geodesic, embedded surface in $M$, whenever $R>0$ is smaller than the minimum of $\tilde{R}$ and the injectivity radius of $M$ at the foot point $\pi_G(P)$ of $P$.
\end{theorem}

\proof Since $R<\tilde{R}$ we may choose $r>0$ such that $R+2r<\tilde{R}$. Now Lemma \ref{lemma82} provides sequences $\ve_i\da 0$ and $q_i\in G_{\ve_i,R,r}(|A_i|^2)$
such that $\lim_{i\to\infty} d_i(q_i, p_i)=0$. Hence we have $\lim_{i\to\infty} E(f_i|B(q_i, R+2r)) = 0$. Since $G_2 M$ is compact there exists a subsequence of $df_i(T_{q_i} F_i)\in G_2 M$ converging to some $P\in G_2 M$. So our claim follows from Proposition \ref{prop91}.
\kasten\vspace{0.2cm}

\begin{remark} \label{rem95}
{\em If we consider the case of $k$-dimensional immersions with $k>2$, our proof will go through under the following additional assumptions:\\
1) A uniform upper bound on ${\mathrm{vol}}_k(B(p_i, \tilde{R}))$.\\
2) An analogue of Proposition \ref{prop63}.}
\end{remark}

The following proposition will be used to prove the existence of {\em complete}, totally geodesic surfaces, see Proposition \ref{prop104}.

\begin{proposition} \label{prop96}
Suppose there are sequences $\ve_i\da 0$, $R_i\to \infty$, and $p_i\in \tilde{G}_{\ve_i, R_i}(|A_i|^2)\subseteq F_i$ such that $\lim_{i\to\infty} E(f_i|B(p_i, R_i)) = 0$ and such that
$\lim_{i\to\infty} df_i(T_{p_i}F_i)=P\in G_2 M$ exists. Given $(v, t)\in SP\times \R$ set $\tilde{P}=\tilde{P}(v, t)={\cal{P}}_t^{c_v}(P)$. Then
$N_{\tilde{P}}(R)= \{c_w(s)|w\in S\tilde{P}, |s|<R\}$ is a totally geodesic, embedded surface in $M$, provided $R>0$ is smaller than the injectivity radius of $M$ at $c_v(t)$.
\end{proposition}

\proof Since $p_i\in \tilde{G}_{\ve_i, R_i}(|A_i|^2)$ we can find sequences $v_i\in S_{p_i} F_i$, $t_i\in \R$ such that $\lim_{i\to\infty} (df_i(v_i), t_i) = (v, t)$ and
$\int_{-R_i}^{R_i} |A_i|^2 \circ c_{v_i}(s)\,ds < \ve_i$ and $c_{v_i}(t_i)\in G_{\ve_i,R_i}(|A_i|^2)$, cf. (\ref{eq37}) and (\ref{eq38}). Now Proposition \ref{prop23}(c) implies that
$\lim_{i\to\infty} df_i(T_{c_{v_i}(t_i)}F_i)= {\cal{P}}_t^{c_v}(P) =\tilde{P}$. Since $c_{v_i}(t_i)\in G_{\ve_i,R_i}(|A_i|^2)=G_{\ve_i,R_i, R_i}(|A_i|^2)$,
we see that the assumptions of Proposition \ref{prop91} are satisfied for the sequence $c_{v_i}(t_i)$ and every choice of $R>0$, $r>0$. So Proposition \ref{prop91} implies our claim.
\kasten

\setcounter{equation}{0}
%%%%%%%%%%%%%%%%%%%%%%%%%%%%%%%%%%%%%%%%%%%%%%%%%%%%%%%%%%%%%%%%%%%%%%%%%%%%%%%
\section{Existence of complete, totally geodesic surfaces} \label{sec10}

We continue to consider a sequence of complete surface immersions $f_i:F_i\to M$ into a compact Riemannian manifold $(M, g)$. %% NEU, S. 10.1
In this section we prove a global version of Proposition \ref{prop91}. Under the assumption that there exist sequences $p_i\in F_i$, $R_i\to \infty$, such that
$\lim_{i\to\infty} E(f_i|B(p_i, R_i)) = 0$, we prove the existence of a complete, totally geodesic surface immersion into $M$. In view of Proposition \ref{prop96}
this amounts to piecing together local totally geodesic surfaces. This is reminiscent of the construction of the leaves of a foliation. Indeed, for general Riemannian manifolds
$(M, g)$, there exists a distribution ${\cal D}_k= {\cal D}_k(M, g)$  in the tangent bundle of the Grassmann bundle $\pi_G: G_k M\to M$ such that the integral manifolds
$L$ of ${\cal D}_k$ correspond to totally geodesic immersions $\pi_G|L: L\to M$. Under our assumptions we will prove that there exists a complete leaf $L\subseteq G_2 M$ of
${\cal D}_2(M, g)$ , i.\,e. a leaf $L$ such that $(\pi_G|L)^{\ast} g$ is a complete Riemannian metric. Then $\pi_G|L: L\to M$ is a complete, totally geodesic surface immersion.\\

We start by collecting some facts concerning integral manifolds of a general distribution ${\cal D}$ of fiber dimension $k$ in the tangent bundle of a general manifold $M$.
We assume that all considered objects are smooth ($=C^\infty$). So, ${\cal D}:M \to G_k M$ is a section of the Grassmann bundle $\pi_G: G_k M\to M$. An {\em integral manifold of }
${\cal D}$ is an immersion $j: L \to M$ of a connected manifold $L$ such that $dj(T_x L)={\cal D}_{j(x)}$ for all $x\in L$. Although we will not use this in our proofs, we note the following important property of integral manifolds.

\begin{remark}\label{rem101} {\em Let $j:L\to M$ be an injective integral manifold of ${\cal D}$. Suppose $P$ is a manifold and $h\in C^\infty(P, M)$ satisfies $h(P)\subseteq j(L)$.
Then $j^{-1} \circ h\in C^\infty(P, L)$}.
\end{remark}

An integral manifold $j:L\to M$ of ${\cal D}$ is called {\em maximal} if $j$ is injective and if the following is true: If $\tilde{j}: \tilde{L}\to M$ is an integral manifold
of ${\cal D}$ and $j(L)\cap \tilde{j}(\tilde{L})\neq \emptyset$, then $\tilde{j}(\tilde{L}) \subseteq j(L)$. If $j: L \to M$ is a maximal integral manifold of ${\cal D}$, then
$j(L)\subseteq M$ is called a {\em leaf of } ${\cal D}$.

One defines a manifold structure on $j(L)$ by declaring $j$ a diffeomorphism. The topology of this manifold structure is finer, and often strictly finer,
than the topology of $j(L)$ as a subspace of $M$. Remark \ref{rem101} shows that this manifold structure of a leaf is in fact independent of the parametrization $j$.

\begin{proposition} \label{prop102}
Suppose $p\in M$ lies in the image of some integral manifold of ${\cal D}$. Then there exists a maximal integral manifold $j:L\to M$ such that
$p\in j(L)$.
\end{proposition}

If ${\cal D}$ is completely integrable, i.\,e. if ${\cal D}$ satisfies the Frobenius condition, then this is proved in textbooks treating foliations.
The more general version stated here admits a similar proof. \\

In particular, if $j_0:L_0\to M$ is an integral manifold of ${\cal D}$, then there exists a leaf $L$ of ${\cal D}$
containing $j_0(L_0)$. Indeed, $q\in L$ iff there exists a finite sequence $j_1:L_1\to M, \ldots, j_n:L_n\to M$ of integral manifolds of ${\cal D}$ such that
$j_{i-1}(L_{i-1}) \cap j_i(L_i)\neq \emptyset$ for $1\leq i\leq n$ and $q\in j_n(L_n)$.\\

Next we briefly recall how $k$-dimensional, totally geodesic immersions into an $m$-dimensional Riemannian manifold $(M, g)$ are related to a $k$-dimensional distribution
${\cal D}_k={\cal D}_k(M, g)$ in the tangent bundle of the Grassmann bundle $G_k M$. So, from now on, $G_k M$ will play the role of the manifold $M$ in the preceding paragraph.
%%%%%%%%%%%%%%%%%%%
Since the bundle $\pi_G: G_k M\to M$ is associated to the principal ${\cal{O}}(m)$-bundle of orthonormal frames there is a natural horizontal distribution ${\cal H}\subseteq T(G_k M)$ induced by the Levi-Civita connection of $g$, cf. \cite{KoNo63}, Chapter II, pp. 87--88. Explicitly ${\cal H}$ is given as follows. If $P\in G_k M$ and $v\in T_{\pi_G(P)}M$, choose a $C^1$-curve
$\gamma:\R\to M$ such that $\dot{\gamma}(0)=v$, and let $P_\gamma(t) = {\cal P}_{t}^\gamma(P)$ denote the parallel transport of $P=P_\gamma(0)$ along $\gamma$.
Then $\dot{P}_\gamma(0)\in T_P(G_k M)$ is independent of the choice of $\gamma$ with $\dot{\gamma}(0) = v$, and defines a linear map
$H_P: T_{\pi_G(P)}M\to T_P(G_k M)$, $H_P(v) = \dot{P}_\gamma(0)$, satisfying $d\pi_G\circ H_P(v)=v$ for all $v\in T_{\pi_G(P)}M$. Then the ($m$-dimensional) horizontal distribution
${\cal H}\subseteq T(G_k M)$ is given by ${\cal H}_P = H_P(T_{\pi_G(P)}M)$ for all $P\in G_k M$. Note that a $C^1$-curve $P:I\to G_k M$ is horizontal, i.\,e. $\dot{P}(t) \in {\cal H}_{P(t)}$ for all $t\in I$, if and only if $P(t)$ is parallel along $(\pi_G\circ P)(t)$. Now we consider the $k$-dimensional distribution ${\cal D}_k={\cal D}_k(M, g)\subseteq {\cal H}$ defined by
$$
{\cal D}_P = H_P(P)
$$
for all $P\in G_k M$.

Note: We have $V\in {\cal D}_P$ if and only if $V\in {\cal H}_P$ and $d\pi_G(V)\in P$.

\begin{lemma}\label{lemma103}
Let $j:N\to M$ be a totally geodesic immersion of a connected, $k$-dimensional manifold $N$, and define $J:N\to G_k M$ by $J(p) = dj(T_pN)$. Then $\pi_G\circ J=j$ and $J$ is an integral manifold of ${\cal D}$, i.\,e. $dJ(T_pN) = {\cal D}_{J(p)}$ for all $p\in N$. Conversely, if $J: N \to G_k M$ is an integral manifold of ${\cal D}$,
then $\pi_G\circ J: N \to M$ is a totally geodesic immersion and $J(p) = d(\pi_G \circ J)(T_p N)$ for all $p\in N$.
\end{lemma}

\proof Suppose first that $j:N \to M$ is a totally geodesic immersion and define $J:N \to G_k M$ by $J(p)=dj(T_p N)$. Let $\tilde{\gamma}$ be a $C^1$-curve in $N$ and $\gamma=j\circ \tilde{\gamma}$. Since $j$ is a totally geodesic immersion we see that $J(\tilde{\gamma}(t)) = dj (T_{\tilde{\gamma}(t)} M)$ is parallel along $\gamma$. By the definition of ${\cal H}_{\gamma(t)}$ this implies that
% $(J\circ \tilde{\gamma})^{\bullet}(t)\in {\cal H}_{J(\gamma(t))}$.
$(J\circ \tilde{\gamma})^{\cdot}(t)\in {\cal H}_{J(\gamma(t))}$. This proves that $dJ(T_p N)\subseteq {\cal H}_{J(p)}$ for all $p\in N$. Additionally we have $d\pi_G \circ dJ = dj$, so that the note preceding Lemma \ref{lemma103} shows that $dJ(T_p N)= {\cal D}_{J(p)}$ for all $p\in N$. Conversely, suppose that $J:N \to G_k M$ satisfies $dJ(T_p N)= {\cal D}_{J(p)}$ for all $p\in N$, and set
$j=\pi_G\circ J$. Then we have for all $p\in N$:
$$
dj(T_p N)= d\pi_G ({\cal D}_{J(p)}) = d\pi_G (H_{J(p)}(J(p))) = J(p).
$$
In particular, $j$ is an immersion. To see that $j$ is totally geodesic note that, for every $C^1$-curve $\tilde{\gamma}$ in N, the curve
$$
t\to dj(T_{\tilde{\gamma}(t)}N) = (J\circ \tilde{\gamma})(t) \in G_k M
$$
is horizontal, i.\,e. parallel along $\gamma=j \circ \tilde{\gamma}$. This proves that $j$ is totally geodesic.\\
\kasten\vspace{0.2cm}

Using Proposition \ref{prop96} and the preceding discussion we will prove:
\begin{proposition} \label{prop104}
Suppose $\ve_i\da 0$, $\rho_i\to \infty$ and $p_i\in \tilde{G}_{\ve_i, \rho_i}(|A_i|^2)$ are sequences such that $\lim_{i\to\infty} E(f_i|B(p_i, \rho_i)) = 0$ and such that
$\lim_{i\to\infty} df_i(T_{p_i}F_i)=P\in G_2 M$ exists. Then there exists a leaf $L$ of ${\cal D}_2(M, g)$ such that $P\in L$ and such that $(\pi_G|L)^\ast g$ is
a complete Riemannian metric on $L$. In particular, $\pi_G|L:L\to M$ is a complete, totally geodesic surface immersion. Moreover, for every $\tilde{P}\in L$ there exists
a sequence $q_i\in G_{\ve_i, \rho_i}(|A|^2)$ such that $\lim_{i\to\infty} df_i(T_{q_i}F_i)=\tilde{P}$ and $d_i(p_i, q_i)\leq d^L(P, \tilde{P})$,
where $d^L$ denotes the distance on $L$ induced by $(\pi_G|L)^\ast g$.
\end{proposition}

\begin{corollary}\label{coroll105}
Suppose $p_i\in F_i$, $R_i\to \infty$ are sequences such that $\lim_{i\to\infty} E(f_i|B(p_i, R_i)) = 0$. Then there exists a complete, totally geodesic surface immersion into $M$.
\end{corollary}

\textsc{Proof of Corollary \ref{coroll105} assuming Proposition \ref{prop104}:}
According to Lemma \ref{lemma83} we can find sequences $\ve_i\da 0$, $\rho_i\to \infty$ and $\rho_i <R_i-1$, and $q_i\in \tilde{G}_{\ve_i, \rho_i}(|A_i|^2)$ such that
$\lim_{i\to\infty} d_i(q_i, {p_i})=0$. Then Proposition \ref{prop104} applies to a subsequence of the sequence $q_i$.
\kasten\vspace{0.2cm}

\textsc{Proof of Proposition \ref{prop104}:} From Proposition \ref{prop96} together with Lemma \ref{lemma103} we conclude the following. For every $(v, t)\in SP\times\R$ there exists an integral manifold of ${\cal D}_2(M, g)$ containing $P_v(t) = {\cal P}_{t}^{c_v}(P)$. Hence, by Proposition \ref{prop102}, there exists a leaf $L$ of ${\cal D}_2(M, g)$ containing $P_v(\R)$.
Moreover, the way $L$ is constructed (resp. Remark \ref{rem101}) implies that for all $v\in SP$ the curve $P_v$ is a smooth curve in $L$. Since $\pi_G\circ P_v$ is the $g$-geodesic $c_v$, we see that $P_v:\R\to L$ is a geodesic with respect to $(\pi_G|L)^\ast g$. Hence every geodesic of $(L, (\pi_G|L)^\ast g)$ through $P$ is defined on all of $\R$. So, by the Hopf-Rinow theorem,
$(L, (\pi_G|L)^\ast g)$ is a complete Riemannian manifold. Lemma \ref{lemma103} implies that $\pi_G|L$ is a totally geodesic immersion. Finally, using the Hopf-Rinow theorem again, we obtain, for every $\tilde{P}\in L$, a geodesic $P_v$, $v\in SP$, such that  $P_v(t)=\tilde{P}$, where $t=d^L(P, \tilde{P})$. Since $p_i\in \tilde{G}_{\ve_i, \rho_i}(|A_i|^2)$, we can find a sequence
$(v_i, t_i)\in S_{p_i} F_i\times \R$ such that $c_{v_i}(t_i)\in {G}_{\ve_i, \rho_i}(|A_i|^2)$, $\lim_{i\to\infty} df_i(v_i)=v$, $t_i\ua t$, and $\int_{-\rho_i}^{\rho_i} |A_i|^2\circ c_{v_i}(s)\,ds<\ve_i$.  Using Proposition \ref{prop23}(c) we see that $\lim_{i\to\infty} df_i(T_{c_{v_i}(t_i)}F_i)=P_v(t) = \tilde{P}$. We set $q_i=c_{v_i}(t_i)$. Then $q_i\in G_{\ve_i, \rho_i}|A|^2)$ and $d_i(p_i, q_i)\leq t_i\leq t= d^L(P, \tilde{P})$, and $\lim_{i\to\infty} df_i(T_{q_i}F_i)= \tilde{P}$.
\kasten\vspace{0.2cm}

If $L$ is a non-compact leaf of ${\cal D}_k(M,g)$ such that $(\pi_G|L)^\ast g$ is complete, one can prove the existence of additional complete leaves in the closure of $L$. In this context the notions 
``lamination'' and ``lamination structure'' from \cite{BaCu17}, Section 2 ($D_1$), seem appropriate.

\begin{proposition} \label{prop106}
Let $(M, g)$ be a compact Riemannian manifold, and suppose $L_0\subseteq G_k M$ is a leaf of ${\cal D}_k(M,g)$ such that $(\pi_G|L)^\ast g$ is complete. Let $S$ denote the closure of $L_0$ in $G_k M$. Then there exists a unique $C^\infty$-lamination structure ${\cal L}$ on $S$ with tangent distribution ${\cal D}_k(M, g)|_S$. If $L$ is a leaf of ${\cal L}$, then $\pi_G|L$ is a
complete, totally geodesic immersion.
\end{proposition}

\begin{remark}\label{rem107}
{\em Since ${\cal D}_k(M, g)|_S$ is the tangent distribution of ${\cal L}$, the leaves of ${\cal L}$ are precisely the leaves of ${\cal D}_k(M, g)$ that are contained in $S$.}
\end{remark}

\begin{remark}\label{rem108}
{\em If $L_0$ is compact, then $L_0=S$ is the only leaf of ${\cal L}$. If all leaves of ${\cal L}$ are noncompact, then ${\cal L}$ will have uncountably many leaves.}
\end{remark}

\textsc{Proof of Proposition \ref{prop106}:} The lamination structure ${\cal L}$ on $S$ will be provided by \cite{BaCu17}, Proposition 2.7, once we know that conditions (a) and (b) in this
proposition hold in our situation. First note that, since $L_0$ is assumed complete, there exists an integral manifold of ${\cal D}_k(M, g)$ through every $P\in S$. This, together with the fact that $S$ is closed, implies conditions (a) and (b). If $L$ is a leaf of ${\cal L}$, then $L$ is a leaf of ${\cal D}_k(M, g)$, see Remark \ref{rem107}, and hence
$\pi_G|L$ is a totally geodesic immersion by Lemma \ref{lemma103}. It is a general fact that the leaves of a lamination (in the sense of \cite{BaCu17}) in a complete Riemannian manifold are complete. In our situation we can also argue that for every leaf $L$ of ${\cal L}$, every $P\in L$ and every $v\in SP$, we have $P_v(t)={\cal P}_t^{c_v}(P)\in L$ for $|t|<\injrad(M, g)$.
This implies that actually $P_v(\R)\subseteq L$. Since the $P_v$, $v\in SP$, are the geodesics of $(L, (\pi_G|L)^\ast g)$ through $P$, we see that $(L, (\pi_G|L)^\ast g)$ is complete.
\kasten\vspace{0.2cm}

\setcounter{equation}{0}
%%%%%%%%%%%%%%%%%%%%%%%%%%%%%%%%%%%%%%%%%%%%%%%%%%%%%%%%%%%%%%%
\section{An upper area estimate for weakly starshaped domains} \label{sec11}

Upper estimates for the area of parallel sets of curves on surfaces under assumptions on the integrated Gaussian curvature go back to work by G.\, Bol \cite{Bol41} and F.\, Fiala \cite{Fia41}.
Rigorous proofs for smooth surfaces were given by P.\,Hartman \cite{Har64}, see also \cite{BuZa88}, Chapter 2, and \cite{SST03}, Chapter 4. We will slightly extend these results from parallel sets to sets that will be called {\em weakly starshaped}. This generalization to weakly starshaped sets is used in the proof of Theorem \ref{thmareabound}. For a proof of Theorem \ref{thrm13}
the known estimate for the area of metric balls is sufficient, see Sections \ref{sec6} and  \ref{sec13}. \\

In this section we will consider a domain $S$ with smooth, compact and connected boundary $\partial S\subseteq S$ on a complete Riemannian surface $(F, g)$. So, the boundary $\partial S$ of $S$ is a simple closed curve that will be denoted by $\Gamma$. We introduce the following notation that will be used throughout this section. The geodesic curvature of $\Gamma$ with respect to the normal pointing into $S$ will be denoted by $\kappa:\Gamma \to \R$. We set $L={\mathrm{length}}(\Gamma)$, ${\cal K}(\Gamma) = \int_\Gamma \kappa(q)\,ds(q)$,
$|{\cal K}|(\Gamma) = \int_\Gamma |\kappa(q)|\,ds(q)$, and ${\cal K}^{+}(\Gamma) = \int_\Gamma \kappa^{+}(q)\,ds(q)$, where $s$ denotes arclength on $\Gamma$.
The distance function $d^{\Gamma}: S\to [0, \infty)$ from $\Gamma=\partial S$ is defined by
$$
d^{\Gamma}(p) = \min\{d(p, x)|x\in \Gamma\},
$$
where $d$ denotes the distance on $F$ induced by $g$. For $t>0$ we set
$$
\Gamma^t = \{p\in S|d^{\Gamma}(p)\leq t\} = (d^{\Gamma})^{-1}([0, t])
$$
and
$$
P \Gamma^t = \{p\in S|d^{\Gamma}(p)= t\} = (d^{\Gamma})^{-1}(\{t\}).
$$
Obviously, we have $\partial{\Gamma}^t \subseteq \Gamma\cup P\Gamma^t$, where strict inclusion may hold for exceptional $t$.

\begin{definition}\label{def111}
A subset $C\subseteq S$ is weakly $\Gamma$-starshaped, if $C$ is closed, $\Gamma\subseteq C$, and the following holds for every $p\in C$ and every geodesic $c:[0, d^{\Gamma}(p)]\to S$ that is a shortest connection from $p$ to $\Gamma$:
$$
c(t) \in \Int(C)\quad\mbox{for all } t\in (0, d^\Gamma(p)).
$$
\end{definition}

\begin{remark}\label{rem112} {\em We do not assume that shortest connections from points in $C$ to $\Gamma$ are unique. So, $\Gamma$ may not be a retract of $C$. This is the reason for the attribute ``weakly''.}
\end{remark}

\begin{remark}\label{rem113} {\em If $C_1$ and $C_2$ are weakly $\Gamma$-starshaped, then so is $C_1\cap C_2$. The collars $\Gamma^t$ of $\Gamma$ are weakly $\Gamma$-starshaped.}
\end{remark}

As in \cite{BuZa88}, Theorem 2.4.2, the area estimate depends on an arbitrary chosen number $k\in\R$. In our application, $k$ will be a negative lower bound for the sectional curvature of the ambient manifold $M$. So we present the estimate only in the case $k<0$, although the cases $k=0$ and $k>0$ can be treated similarly, see \cite{BuZa88}, Theorem 2.4.2, for the corresponding formulae. If $B\subseteq F$ is measurable we set
$$
\omega_k^{-}(B) =\int_B(K-k)^{-}d\volz,
$$
where $K$ denotes the Gaussian curvature of $F$. For $C\subseteq S$ and $t>0$ we set
$$
C^t = C\cap \Gamma^t \mbox{ and } PC^t=C\cap P\Gamma^t = \{p\in C| d^\Gamma(p)=t\}.\\
$$

The aim of this section is to prove:
\begin{proposition} \label{prop114}
Suppose $t>0$, $k<0$, and $C\subseteq S$ is weakly $\Gamma$-starshaped. Then
$$
\volz(C^t) \leq \frac{1}{-k}({\cal K}^{+}(\Gamma)+ \omega_k^{-}(C^t))(\cosh(\sqrt{-k}t)-1) + \frac{L}{\sqrt{-k}}\sinh(\sqrt{-k}t)
$$
and
$$
{\cal H}^1(PC^t) \leq \frac{1}{\sqrt{-k}}({\cal K}^{+}(\Gamma) + \omega_k^{-}(C^t))\sinh(\sqrt{-k}t) +L\cosh(\sqrt{-k}t).
$$
\end{proposition}

Note: In the case $C=\Gamma^t$ and $\kappa\geq 0$, the first estimate reduces to \cite{BuZa88}, Theorem 2.4.2, case $k<0$.
      If $ \kappa\geq 0$ is not assumed, it is weaker, since the estimate in \cite{BuZa88}, Theorem 2.4.2, contains the term ${\cal K}(\Gamma)$ instead of ${\cal K}^{+}(\Gamma)$.\\

Specializing Definition \ref{def111} we say that a closed subset $C\subseteq F$ is {\em weakly } $p$-{\em starshaped}, if $p\in \Int(C)$ and every shortest geodesic
$c:[0, d(p, q)]\to F$ from a point $q=c(0)\in C$ to $p$ satisfies $c(t)\in\Int(C)$ for all $t\in(0, d(p,q)]$. Then $C$ is weakly $\Gamma$-starshaped, where
$S=F\backslash B(p, \ve)$, $\Gamma = \partial S=\partial B(p, \ve)$, provided $\overline{B(p,\ve)}\subseteq \Int(C)$ and $\ve$ is smaller than the injectivity radius at $p$. So, in the limit
$\ve\da 0$ Proposition \ref{prop114} implies:

\begin{corollary}\label{coroll115} If $C\subseteq F$ is weakly $p$-starshaped, and $C\subseteq B(p,r)$ for some $r>0$, then the following holds for every $k<0$:
$$
\volz(C)\leq \frac{1}{-k}(2\pi+ \omega_k^{-}(C))(\cosh(\sqrt{-k}r)-1).
$$
\end{corollary}

Next we assume that $f:F\to M$ is a complete surface immersion into a Riemannian manifold with sectional curvature bounded below by $k<0$. Then the Gau{\ss} equation implies
$K\geq k-\frac{1}{2}|A|^2$, cf. (\ref{eqgauss}), and hence
\begin{equation}\label{eq111}
\omega_k^{-}(B) \leq E(f|B)
\end{equation}
for every measurable subset $B\subseteq F$. So, in this situation, the inequality in Corollary \ref{coroll115} says:
\begin{equation}\label{eq112}
\volz(C)\leq \frac{1}{-k}(2\pi+ E(f|C))(\cosh(\sqrt{-k}r)-1).
\end{equation}
Finally, we note that $C=\overline{B(p,r)}$ is weakly $p$-starshaped and $\partial B(p, r)\subseteq PC^r$. Hence Proposition \ref{prop114} together with (\ref{eq111}) imply:
\begin{equation}\label{eq113}
{\cal H}^1(\partial B(p, r)) \leq \frac{1}{\sqrt{-k}}(2\pi+ E(f|B(p, r)))\sinh(\sqrt{-k}r).
\end{equation}

At first sight one may have the impression that the proof given in \cite{Har64} covers the more general case of a weakly $\Gamma$-starshaped set $C$. However, there are some technical problems
concerning the analytical properties of the function $t\to {\cal H}^1(PC^t)$. To circumvent these problems we first approximate $C$ by a weakly $\Gamma$-starshaped subset of $C$ that does not intersect the cut locus of $\Gamma$, and that has a smooth and generic boundary. For such sets the original ideas from  \cite{Bol41} and \cite{Fia41} apply directly, and Proposition \ref{prop114} will follow from this by approximation.\\

In order to describe this approximation process we introduce the following tools that are taken from \cite{SST03}, Chapter 4. We let $N$ denote the unit normal field along $\Gamma$ that points into $S$ and define
$$
Z:\Gamma\times[0, \infty)\to F, \quad Z(p,t)=\exp(t N(p)).
$$
The distance function $\rho$ to the cut locus of $\Gamma$,
$$
\rho:\Gamma \to (0, \infty], \quad \rho(p) =\sup\{t|d^\Gamma(Z(p,t))=t\}
$$
is continuous. The set
$$
CL(\Gamma)=\{Z(p, \rho(p))|p\in\Gamma, \rho(p)<\infty\}
$$
is the cut locus of $\Gamma$ (in $S$). If we restrict $Z$ to the set $D=\{(p,t)|p\in\Gamma, 0\leq t <\rho(p)\}$, then $Z|D$ is a diffeomorphism onto $S\backslash CL(\Gamma)$.\\

It obviously suffices to prove Proposition \ref{prop114} for {\em compact} weakly $\Gamma$-starshaped sets $C$, since we can replace $C$ by $C\cap\Gamma^r$ for $r>t$.
So, in the sequel we will assume that $C\subseteq \Gamma^r$ for some $r>0$. We then define a function $g=g_C:\Gamma\to (0, r]$ by
$$
g(p)=\begin{cases}
\rho(p)& \mbox{if } \rho(p) < \infty \mbox{ and } Z(p, \rho(p))\in C\\
\inf\{t>0|Z(p,t)\notin C\}& \mbox{otherwise.}
\end{cases}
$$
\begin{remark}\label{rem116}
{\em Since $C$ is weakly $\Gamma$-starshaped, $Z$ maps the set $\{(p,t)|p\in \Gamma, 0<t<g(p)\}\subseteq D$ diffeomorphically onto $\Int(C)\backslash CL(\Gamma)$, while
$Z(\graph(g))=(\partial C\backslash\Gamma)\cup(C\cap CL(\Gamma))$.}
\end{remark}
In particular, Remark \ref{rem116} implies
\begin{equation}\label{eq114}
\volz(\partial C)=0.
\end{equation}
\begin{lemma}\label{lemma117}
The function $g_C:\Gamma\to (0, r]$ is continuous.
\end{lemma}

\proof Suppose the sequence $p_i\in \Gamma$ converges to $p\in \Gamma$. Since $C$ is closed and $Z(p_i, g(s_i))\in C$ for all $i\in \N$, we have $Z(p, t)\in C$ for every limit point $t$ of the sequence $g(p_i)\in (0, r]$. Since $g(p_i) \leq \rho(p_i)$ we conclude that $t\leq \lim_{i\to\infty}\rho(p_i)=\rho(p)$. Hence the definition of $g$ implies $t\leq g(p)$. So, to prove continuity of $g$, it suffices to show that the assumption $t< g(p)$ leads to a contradiction. Indeed, if $t<\tau<g(p)$ then $Z(p,\tau)\in \Int(C)$, and hence $Z(p_i, \tau)\in C$ for almost all $i\in \N$.
Since $\tau < g(p) \leq \rho(p) = \lim_{i\to\infty}\rho(p_i)$, this implies $\tau\leq g(p_i)$ for almost all $i\in \N$. Since $t$ is a limit point of the $g(p_i)$, this contradicts our assumption $t<\tau$.
\kasten\vspace{0.2cm}

Given $t>0$ we let $D_1 Z(p,t)$ denote the derivative of the curve $p\in\Gamma\to Z(p,t)$ with respect to arclength on $\Gamma$. Then, if $A\subseteq \Gamma$ is measurable, we have
\begin{equation}\label{eq115}
{\cal H}^1 (Z(A\times\{t\})) \leq \int_A |D_1 Z(p,t)|\,ds(p),
\end{equation}
cf. \cite{Fed69}, Corollary 3.2.20. Using Remark \ref{rem116} we see that $PC^t = Z(A_C^t\times\{t\})$, where $A_C^t=\{p\in \Gamma|g_C(p)\geq t\}$. We set
$$
L_C(t) = \int_{A_C^t}|D_1 Z(p,t)|\,ds(p).
$$
From (\ref{eq115}) we obtain
\begin{equation}\label{eq116}
{\cal H}^1 (PC^t) \leq L_C(t).
\end{equation}

\begin{lemma}\label{lemma118} Let $C\subseteq \Gamma^r$ be weakly $\Gamma$-starshaped, and $\delta>0$. Then there exists a weakly $\Gamma$-starshaped set
$\tilde{C}\subseteq C\backslash CL(\Gamma)$ with the following properties:\\
\begin{itemize}
\item[\em{(i)}]   $g_C(p)-\delta < g_{\tilde{C}}(p)$ for all $p\in \Gamma$.
\item[\em{(ii)}]  $g_{\tilde{C}}:\Gamma \to (0, r)$ is a smooth Morse function.
\item[\em{(iii)}] $\volz (C\backslash \tilde{C})\leq l^2 L\delta$, where $l$ denotes a common Lipschitz constant for $Z|\Gamma\times [0, r]$ and $D_1 Z|\Gamma\times [0, r]$.
\item[\em{(iv)}]  $L_C(t) \leq L_{\tilde{C}}(t-\delta)+l L\delta$, whenever $t>\delta$.
\end{itemize}
\end{lemma}

\proof
Since, by Lemma \ref{lemma117}, $g_C$ is continuous, it is a standard result from Morse theory that there exists a smooth Morse function $\tilde{g}:\Gamma\to(0, r)$ such that
$g_C(p)-\delta <\tilde{g}(p) < g_C(p)$ for all $p\in \Gamma$. Then we set
$$
\tilde{C} = \{Z(p,t)|p\in \Gamma, 0\leq t\leq \tilde{g}(p)\}.
$$
Since $\tilde{g}< g_C\leq \rho$ we have $\tilde{C}\subseteq C\backslash CL(\Gamma)$. Moreover, $\tilde{C}$ is weakly $\Gamma$-starshaped and $g_{\tilde{C}}=\tilde{g}$.
This takes care of properties (i) and (ii). To prove (iii) note that $C\backslash \tilde{C}= \{Z(p,t)|\tilde{g}(p)< t \leq g_C(p)\}$, cf. Remark \ref{rem116}.
Since $Z|\Gamma \times [0, r]$ is $L$-Lipschitz we obtain (iii). Finally note that $g_C -\delta < \tilde{g}$ implies that $A_C^t \subseteq A_{\tilde{C}}^{t-\delta}$ for all $t>\delta$.
Hence, if $t>\delta$ then
$$
\begin{array}{rcl}
\D L_C(t)=\int_{A_C^t} |D_1 Z(p,t)|\,ds(p)&\leq&\\
\D \int_{A_{\tilde{C}}^{t-\delta}}|D_1 Z(p,t-\delta)|\,ds(p)+\int_\Gamma |D_1 Z(p,t)-D_1 Z(p,t-\delta)|\,ds(p)&\leq& L_{\tilde{C}}(t-\delta)+lL\delta.
\end{array}\vspace{-0.2cm}
$$
\kasten

Note that $P\Gamma^t\backslash CL(\Gamma)$ is a smooth 1-dimensional submanifold of $F$, if $P\Gamma^t \backslash CL(\Gamma)\neq \emptyset$.
We let $\kappa^t: P\Gamma^t \backslash CL(\Gamma)\to\R$ denote the geodesic curvature of $P\Gamma^t \backslash CL(\Gamma)$ with respect to the outward pointing unit normal
$\grad( d^\Gamma)$.

\begin{lemma}\label{lemma119} Let $C\subseteq \Gamma^r \backslash CL(\Gamma)$ be weakly $\Gamma$-starshaped, and assume that $g=g_C$ is a smooth Morse function. Then the function
$v:(0, \infty)\to (0, \infty)$, $v(t) = \volz(C^t)$ is $C^1$ with derivative $v'= L_C$, and $v$ is smooth on the set of regular values of $g$. If $t>0$ is a regular value of $g$,
then $v''(t)\leq \int_{PC^t} \kappa^t(q) \,ds(q)$.
\end{lemma}

\proof We fix a maximal interval $(t_{-}, t_{+})\subseteq (0, \infty)$ of regular values of $g$. First we prove that $L_C|(t_{-}, t_{+})$ is smooth. For every $t\in(t_{-}, t_{+})$
the set $A_C^t=g^{-1}([t,\infty))\subseteq \Gamma$ is a finite union of disjoint intervals $I_i(t)$, $1\leq i \leq l$, and
\begin{equation}\label{eq117}
Z(\bigcup_{i=1}^l (I_i(t)\times \{t\})) = PC^t.
\end{equation}
Moreover, $l$ is independent of $t\in (t_{-}, t_{+})$, and the endpoints of the intervals $I_i(t)\subseteq \Gamma$ depend smoothly on $t\in (t_{-}, t_{+})$.
Since $C\cap CL (\Gamma)=\emptyset$, the map $Z$ is diffeomorphic when restricted to the set $\{(p, t)|p\in \Gamma$, $0\leq t\leq g(p)\}$. In particular, for
$t\in (t_{-}, t_{+})$, $Z$ maps $I_i(t)\times\{t\}$, $1\leq i \leq l$, diffeomorphically onto disjoint $\mathrm{arcs\,}PC_i^t\subseteq \Gamma^t \backslash CL(\Gamma)$ and
$\bigcup_{i=1}^l PC_i^t= PC^t$ by (\ref{eq117}). So, for $t\in (t_{-}, t_{+})$, we have
$$
L_C(t) = \sum_{i=1}^l {\mathrm{length}}(PC_i^t) =  \sum_{i=1}^l \int_{I_i(t)} |D_1 Z(p, t)|\,ds(p).
$$
Since $D_1 Z(p, t)\neq 0$ for $p\in I_i(t)$, we conclude that $L_C|(t_{-}, t_{+})$ is smooth. Next we estimate $L_C'$. Since $h\to Z(p, t+h)$ parametrizes geodesics orthogonal to $PC^t$,
standard differential geometry shows that
\begin{equation}\label{eq118}
\lim_{h\to 0}\frac{1}{h} \int_{I_i(t)}|D_1 Z(p, t+h)|-|D_1 Z(p, t)|\,ds(p)=\int_{PC_i^t} \kappa^t(q)\,ds(q).
\end{equation}
Now suppose $h>0$ and $[t, t+h]\subseteq (t_{-}, t_{+})$. Then $A_C^{t+h} = g^{-1}([t+h, \infty))\subseteq g^{-1}([t, \infty)) = A_C^t$, and hence $I_i(t+h) \subseteq I_i(t)$ for
$1\leq i \leq l$. Using (\ref{eq118}) we obtain
\begin{eqnarray}
\label{eq119} L_C'(t) &=& \lim_{h\da 0}\frac{1}{h} (\sum_{i=1}^l (\int_{I_i(t+h)}|D_1 Z(p, t+h)|\,ds(p)- \int_{I_i(t)}|D_1 Z(p, t)|\,ds(p)))\\
\nonumber             &\leq& \int_{PC^t} \kappa^t(q)\,ds(q).
\end{eqnarray}
Next we show that $v$ is $C^1$, and $v' = L_C$. Since $g''(p) \neq 0$ at every singular point $p\in \Gamma$ of $g$, one easily concludes that $L_C$ is continuous on all of $(0,\infty)$.
Moreover, $PC^t = C\cap P\Gamma^t$ are the level sets of the distance function $d^\Gamma$ restricted to $C$. Hence the coarea formula, see e.\,g. \cite{Fed69}, Theorem 3.2.22, implies that $v(t) = \int_0^t L_C(\tau)\,d\tau$. Since $L_C$ is continuous we see that $v$ is $C^1$ and $v' = L_C$. Since $L_C$ is smooth on the set of regular values of $g$, the same holds for $v$.
Moreover, (\ref{eq119}) implies that $v''(t) \leq \int_{PC^t} \kappa^t(q)\,ds(q)$ for every regular value $t$ of $g$.
\kasten\vspace{0.2cm}

\textsc{Proof of Proposition \ref{prop114}:} First we treat the case that $C$ satisfies the assumptions in Lemma \ref{lemma119}. Then the general case will easily follow from Lemma \ref{lemma118}. We intend to use the Gau{\ss}-Bonnet formula to prove that the following differential inequality for the function $v(t) = \volz(C^t)$ is valid for regular values $t$ of $g=g_C$.
\begin{equation}\label{eq1110}
v''(t) + kv(t) \leq {\cal K}^{+}(\Gamma)+\omega_k^{-}(C^t).
\end{equation}
As in the proof of Lemma \ref{lemma119} we consider a maximal interval $(t_{-}, t_{+})$ of regular values of $g$, and find $l\in \N$ and disjoint intervals $I_i(t)$, $1\leq i\leq l$,
in $\Gamma$ such that $g^{-1}([t, \infty))=\bigcup_{i=1}^l I_i(t)$. Then the sets $C_i^t = Z(I_i(t)\times [0, t])$, $1\leq i\leq l$, are disjoint subsets of $C^t$. Each $C_i^t$ is a simply connected domain with piecewise smooth boundary consisting of $PC_i^t = Z(I_i(t)\times \{t\})$, $I_i(t)$, and the two geodesic arcs $Z(\partial I_i(t)\times [0, t])$.
Since $\partial C_i^t$ has right angles at its four corners, the Gau{\ss}-Bonnet formula applied to $C_i^t$ and summation over $i$ imply
$$
\int_{\bigcup_{i=1}^l C_i(t)} K d\volz + \int_{PC^t} \kappa^t(q)\,ds(q) - \int_{\bigcup_{i=1}^l I_i(t)} \kappa(p)\,ds(p)=0.
$$
Now the preceding equality, Lemma \ref{lemma119}, and the obvious inequality
$$
kv(t)\leq \int_{\bigcup_{i=1}^l C_i(t)} K d\volz +\omega_k^{-}(C^t)
$$
combined prove (\ref{eq1110}).\\

Next we show that the differential inequality (\ref{eq1110}) implies the estimates claimed in Proposition \ref{prop114}.
We fix $t>0$. The function $v|[0, t]$ is $C^1$ and smooth except at the finitely many critical points of the Morse function $g_C$.
For every $a\geq {\cal K}^{+}(\Gamma)+\omega_k^{-}(C^t)$ and every regular value $\tau\in[0, t]$ of $g_C$, $v$ satisfies the differential inequality
$$
v''(\tau)+ kv(\tau) \leq a,
$$
cf. (\ref{eq1110}). This allows us to compare $v|[0, t]$ and $v'|[0, t]= L_C|[0, t]$ to the solution
$w(\tau)=w_{k, a, L}(\tau) = \frac{a}{-k}(\cosh(\sqrt{-k}\tau)-1)+ \frac{L}{\sqrt{-k}} \sinh(\sqrt{-k}\tau)$ of the initial value problem $w'' + kw = a$, $w(0)=0$, $w'(0)=L$.
So, if $a\geq {\cal K}^{+}(\Gamma)+\omega_k^{-}(C^t)$, we obtain
\begin{equation}\label{eq1111}
v(\tau)\leq w_{k, a, L}(\tau) \mbox{ and } v'(\tau)\leq w'_{k, a, L}(\tau)\quad \mbox{for all } \tau\in[0, t].
\end{equation}

Finally, if $C$ is only weakly $\Gamma$-starshaped and $t>0$ and $\delta\in(0, t)$, we use Lemma \ref{lemma118} to approximate
$C^t$ by a set $\tilde{C}=\tilde{C}_\delta \subseteq C^t$ that satisfies the assumptions of Lemma \ref{lemma119}. We set $a= {\cal K}^{+}(\Gamma)+\omega_k^{-}(C^t)$. Then
$a\geq {\cal K}^{+}(\Gamma)+\omega_k^{-}(\tilde{C})$ and $\tilde{C}=\tilde{C}^{t}$, so that (\ref{eq1111}) implies $\volz(\tilde{C})\leq w_{k, a, L}(t)$ and
$L_{\tilde{C}}(\tau)\leq w'_{k, a, L}(\tau)$ for $\tau\in[0, t]$. Letting $\delta\da 0$ and using Lemma \ref{lemma118} (iii) and (iv), we see that
$\volz(C^t)\leq w_{k, a, L}(t)$ and $L_C(t) \leq w'_{k, a, L}(t)$. Since ${\cal H}^1(PC^t) \leq L_C(t)$ by (\ref{eq116}), the preceding inequalities prove Proposition \ref{prop114}.
\kasten

\setcounter{equation}{0}
%%%%%%%%%%%%%%%%%%%%%%%%%%%%%%%%%%%%%%%%%%%%%%%%%%%%%%%%%%%%%%%
\section{Proof of Theorem \ref{thmareabound}} \label{sec12}

In this section we prove Theorem \ref{thmareabound} by combining Corollary \ref{coroll105} with Corollary \ref{coroll123} below.
The main new result is the following Proposition \ref{prop121} that is based on the results of the preceding section.\\

For $k<0$ we abbreviate $\frac{-1}{k}(\cosh(\sqrt{-k}r)-1)$ by $g_k(r)$.
\begin{proposition} \label{prop121} Let $f: F \to M$ be an isometric immersion of a compact Riemannian surface $(F, g)$ into a Riemannian manifold $(M, \bar{g})$
with sectional curvature bounded below by $k<0$. If $R>0$ and $E(f) g_k(R) < \volz(F)$, then there exists $p\in F$ such that
$$
E(f|B(p, R/2)) \leq 2\pi\frac{E(f) g_k(R)}{\volz(F) - E(f) g_k(R)}.
$$
\end{proposition}

\begin{remark}\label{rem122}
{\em If $(F, g)$ is noncompact and complete, and if $E(f) < \infty$, then $\inf_{p\in F} E(f|B(p, R))=0$ holds for every $R>0$.}
\end{remark}

\begin{corollary}\label{coroll123}
Let $f_i: F_i \to M$ be a sequence of complete surface immersions into a Riemannian manifold $(M, \bar{g})$ with sectional curvature bounded below. If $E(f_i)<\infty$ for all $i\in \N$
and $\lim_{i\to \infty}E(f_i)/\volz(F_i)=0$, then there exist sequences $p_i\in F_i$ and $R_i\to \infty$ such that $\lim_{i\to \infty}E(f_i|B(p_i, R_i))=0$.
\end{corollary}

\textsc{Proof of Corollary \ref{coroll123} assuming Proposition \ref{prop121}:}\\
Since $\lim_{i\to \infty}E(f_i)/\volz(F_i)=0$ we can find a sequence $R_i\to \infty$ such that
$$
\lim_{i\to \infty} g_k(2R_i)E(f_i)/\volz(F_i)=0,
$$
where $k<0$ denotes a lower bound for the sectional curvature of $(M, \bar{g})$. Then Proposition \ref{prop121} together with Remark \ref{rem122} provide points $p_i\in F_i$ such that
$\lim_{i\to \infty}E(f_i|B(p_i, R_i))=0$.
\kasten

To prove Proposition \ref{prop121} we use a decomposition of $F$ into Voronoi cells stemming from an $R$-net on $F$, and apply Proposition \ref{prop114} to the Voronoi cells.
First we recall some facts concerning these concepts. Given $R>0$ we consider an $R$-net $N\subseteq F$, i.\,e. we have $\bigcup_{p\in N} B(p, R)=F$ and $d(p, q)\geq R$ for all $p\neq q$ in $N$. Then the Voronoi cell of $p\in N$ is defined by
$$
C_p = C_p(N)=\{x\in F|d(p, x)=\min_{q\in N} d(q, x)\}.
$$
As direct consequences of these definitions we have:
\begin{equation}\label{eq121}
\bigcup_{p\in N} C_p= F
\end{equation}
\begin{equation}\label{eq122}
\mbox{If } p\in N \mbox{ then } \overline{B(p, R/2)} \subseteq C_p\subseteq B(p, R).
\end{equation}

\begin{lemma}\label{lemma124}
If $p\in N$ then $C_p$ is weakly $p$-starshaped, and $\Int(C_p) = \{x\in F|d(p, x) < d(q, x)$ for all $q\in N\backslash \{p\}\}$.
\end{lemma}

The proof of Lemma \ref{lemma124} relies on the following fact that is a direct consequence of the regularity of arclength-parametrized shortest connections.\\

{\bf Fact 12.5 } Let $c:[0, r]\to F$ be a shortest geodesic from $c(0)=x$ to $c(r) =p$, $r=d(p, x)$. If $q\in F\backslash \{p\}$ and $d(q,x)=d(p, x)$, then
$d(p, c(s))< d(q, c(s))$ for all $s\in (0, r]$.\\

\proof First, we have $d(p, c(s))= r-s = d(q, x)-s\leq d(q, c(s))$ by the triangle inequality. Now, contrary to our claim, assume that
$d(p, c(s))= d(q, c(s))$, and let $\bar{c}:[0, r-s]\to F$ be a shortest geodesic from $\bar{c}(0) = c(s)$ to $\bar{c}(r-s)=q$. Joining $\bar{c}$ to $c|[0, s]$ we obtain an
arclength-parametrized curve of length $r= d(q, x)$ from $x$ to $q$. This implies that this curve is a geodesic, in particular
$\dot{\bar{c}}(0) = \dot{c}(s)$ and $c(t) = {\bar{c}}(t-s)$ for all $t\in [0, r]$. So, $p=c(r) = \bar{c}(r-s)=q$, in contradiction to our assumption $p\neq q$.
\kasten\vspace{0.2cm}

\textsc{Proof of Lemma \ref{lemma124}:} We start by proving the statement concerning $\Int(C_p)$. If $x\in F$ and $d(p, x) < d (q, x)$ for all $q\in N\backslash \{p\}$,
then there exists $\ve >0$ such that $d(p, x) < d(q,x)-\ve$ for all $q\in N\backslash \{p\}$, since $N$ is discrete. This implies $B(x, \ve)\subseteq C_p$, so that $x\in \Int(C_p)$.
Conversely, assume that $x\in \Int(C_p)$ and $d(p, x) = d(q, x)$ for some $q\in N\backslash \{p\}$. Reversing the role of $p$ and $q$ in Fact 12.5 we let $c$ be a shortest geodesic from
$x=c(0)$ to $q$. Then Fact 12.5 implies that $d(q, c(s)) < d(p, c(s))$ for $s\in (0, d(x, q)]$. In particular, we have $c(s)\notin C_p$ for small $s>0$, in contradiction to
$c(0) = x\in \Int(C_p)$. Finally, we prove that $C_p$ is weakly $p$-starshaped. Suppose $x\in C_p$, $r= d(x, p)$, and $c:[0, r]\to F$ is a shortest geodesic from $x=c(0)$ to $p=c(r)$.
We intend to show that, for every $q\in N\backslash \{p\}$, we have $d(c(s), p)< d(c(s), q)$ for all $s\in (0, r]$. This will prove $c(s)\in \Int(C_p)$ for all $s\in (0, r]$ by the preceding argument. If $q\in N\backslash \{p\}$ and $r=d(p, x) < d(q, x)$, then $d(p, c(s)) = r-s < d(q, x)-s \leq d(q, c(s))$ for all $s\in (0, r]$. If $q\in N\backslash \{p\}$ and $d(p, x)= d(q, x)$,
then Fact 12.5 implies $d(c(s), p) < d(c(s), q)$ for all $s\in (0, r]$.
\kasten\vspace{0.2cm}

As a consequence of the characterization of $\Int(C_p)$ in Lemma \ref{lemma124} we obtain
\begin{equation}\label{eq123}
\mbox{If } p\neq q \mbox{ are in } N \mbox{ then } C_p\cap C_q = (\partial C_p)\cap (\partial C_q).
\end{equation}

We note that (\ref{eq121})-(\ref{eq122}), Lemma \ref{lemma124}, and Fact 12.5, are true in more general contexts, e.\,g. in complete, symmetric Finsler manifolds of arbitrary dimension.\\

From (\ref{eq114}) and (\ref{eq123}) we obtain $\volz(C_p\cap C_q) = 0$ for all $p\neq q$ in $N$, and hence (\ref{eq121}) implies
\begin{equation}\label{eq124}
\volz(F) = \sum_{p\in N} \volz(C_p).
\end{equation}

\textsc{Proof of Proposition \ref{prop121}:} We set $\delta = \min_{p\in N} E(f|C_p)/\volz(C_p)$. Then we use (\ref{eq121}) and (\ref{eq124}) to see that $E(f) \geq \delta \volz(F)$,
i.\,e.
\begin{equation}\label{eq125}
\delta \leq E(f)/\volz(F).
\end{equation}
In particular, our assumption implies $\delta g_k(R) < 1$. Choose $p\in N$ such that $E(f|C_p) = \delta \volz(C_p)$. Using Lemma \ref{lemma124},  (\ref{eq112}) and (\ref{eq122}),
we obtain
$$
E(f|C_p) = \delta \volz(C_p) \leq \delta(2\pi + E(f|C_p))g_k(R)
$$
and, equivalently,
$$
E(f|C_p) (1-\delta g_k(R)) \leq 2\pi \delta g_k(R).
$$
According to (\ref{eq122}) and (\ref{eq125}) this implies
$$
E(f|B(p, R/2))\leq 2\pi\frac{E(f) g_k(R)}{\volz(F) - E(f) g_k(R)}.
$$
\kasten\vspace{0.2cm}

\textsc{Proof of Theorem \ref{thmareabound}:} We argue by contradiction, and assume that there exists a sequence $f_i:F_i\to M$ of complete surface immersions into a compact Riemannian
manifold $M$ such that $E(f_i)<\infty$ for all $i\in \N$ and $\lim_{i\to \infty}E(f_i)/\volz(F_i) = 0$. Then Corollary \ref{coroll123} provides sequences $p_i\in F_i$ and $R_i\to \infty$
such that $\lim_{i\to \infty}E(f_i|B(p_i, R_i))=0$. Now we can use Corollary \ref{coroll105} to obtain a complete, totally geodesic surface immersion into $M$, in contradiction to the
assumption made in Theorem \ref{thmareabound}.
\kasten

\setcounter{equation}{0}
%%%%%%%%%%%%%%%%%%%%%%%%%%%%%%%%%%%%%%%%%%%%%%%%%%%%%%%%%%%%%%%
\section{Hausdorff convergence and proof of Theorem \ref{thrm13}} \label{sec13}

We continue to consider a sequence of complete surface immersions $f_i:F_i\to M$ into a compact Riemannian manifold $(M, g)$. We will assume that $p_i \in F_i$ is a sequence satisfying the assumptions of Proposition \ref{prop104}, i.\,e. there exist sequences $\ve_i\da 0$, $\rho_i\to \infty$ such that
\begin{align}
\label{eq131} \D \lim_{i\to \infty}E(f_i|B(p_i, \rho_i))=0, \text{ and}\\
\label{eq132} \D p_i\in \tilde{G}_{\ve_i, \rho_i}(|A_i|^2), \text{ and}\\
\label{eq133} \D \lim_{i\to \infty}df_i(T_{p_i}F_i) = P_0\in G_2M.
\end{align}

Then Proposition \ref{prop104} provides a leaf $L$ of ${\cal D}_2(M, g)$ such that $P_0\in L$ and $\pi_G|L:L\to M$ is a complete, totally geodesic surface immersion.
Under these assumptions we will prove
\begin{theorem} \label{thrm131}
For every $R>0$ the sequence of compact sets $f_i(\overline{B(p_i, R)})\subseteq M$ Hausdorff converges to $\overline{N_{P_0}(R)}=\{c_v(t) |v\in SP_0, |t|\leq R\}$.
\end{theorem}

\begin{remark}\label{rem132}
{\em If $\overline{L_{P_0}(R)}$ denotes the closed metric ball in $(L, (\pi_G|L)^\ast g)$ with center $P_0\in L$ and radius $R>0$, then
$\pi_G(\overline{L_{P_0}(R)}) = \overline{N_{P_0}(R)}$.}
\end{remark}

For the proof of Theorem \ref{thrm13} we rely on the results on pointed Gromov-Hausdorff convergence from Section \ref{sec7}, in particular on Propositions \ref{prop73}(b) and
\ref{prop75}(b). So, we will choose a subsequence, denoted by the same symbols, such that the sequence of pointed metric spaces $(F_i, p_i)_{i\in \N}$ converges to a proper length space
$(Y, y_0)$ with respect to (definite) Gromov-Hausdorff convergence, and such that the immersions $f_i$ converge to a 1-Lipschitz map $f:(Y, y_0)\to (M, p_0)$ where $p_0 = \pi_G(P_0)$.
The distance functions on $F_i$ will be denoted by $d_i$, and the distance functions on $Y$ resp. $M$ by $d$ resp. $d^M$. The distance functions on $F_i\cup Y$ determining the definite convergence 
will be denoted by $\delta_i$.\\

The following proposition will play a crucial role in the proof of Theorem \ref{thrm13}. Its proof is given following Corollary \ref{coroll1312}.
\begin{proposition} \label{prop133}
There exists a covering map $\tilde{f}:(Y, y_0)\to (L, P_0)$ such that $\pi_G\circ \tilde{f}=f$. In fact, $\tilde{f}$ is a local isometry from $(Y, d)$ onto $L$ with the distance induced by $(\pi_G|L)^\ast g$.
\end{proposition}

In this section we will call a sequence $q_i\in F_i$ {\em good} resp. {\em{very good} } if there exist sequences $\ve_i\da 0$, $\rho_i\to \infty$ and such that
$q_i \in G_{\ve_i, \rho_i}(|A_i|^2)$ resp. $q_i \in \tilde{G}_{\ve_i, \rho_i}(|A_i|^2)$, see Definition \ref{def35}, (\ref{eq37}) and (\ref{eq38}).
As a consequence of Lemma \ref{lemma83} we have:
\begin{equation}\label{eq134}
\text{Every } y\in Y \text{ is the limit of a very good sequence } q_i\in F_i.
\end{equation}

In a series of lemmas we will analyse the properties of the limit map $f:(Y, y_0) \to (M, p_0)$. This will lead to a proof of Proposition \ref{prop133} and Theorem \ref{thrm131}.

\begin{lemma}\label{lemma134}
Suppose $q_i\in F_i$ is a good sequence converging to $y\in Y$, and $P\in G_2M$ is a limit plane of the sequence $df_i(T_{q_i}F_i)$. If $0< r < \injrad(M)$, then
$f(\partial B(y, r)) \supseteq S_P(r)=\{c_v(r)|v\in SP\}$. Moreover, $f(B(y, r)) \supseteq N_P(r)=\{c_v(t)|(v, t)\in SP\times[0, r)\}$ for all $r>0$.
\end{lemma}

\proof Since $P$ is a limit plane of the sequence $df_i(T_{q_i} F_i)$, we can find a subsequence, denoted by the same symbols, and, for every $v\in SP$, a sequence
$w_i\in S_{q_i} F_i$ such that $\lim_{i\to \infty} df_i(w_i)=v$. Since $q_i$ is a good sequence we can approximate the $w_i$ by $v_i\in S_{q_i} F_i$ such that
$\lim_{i\to \infty} df_i(v_i) =\lim_{i\to \infty} df_i(w_i) =v$ and $\lim_{i\to \infty} \int_{-R}^R |A_i|^2 \circ c_{v_i}(t)\,dt= 0$ for every $R>0$.
We choose $R$ such that $r< R < \injrad (M)$. Proposition \ref{prop23} implies that the curves $f_i\circ c_{v_i} |[-R, R]$ \, $C^1$-converge to $c_v|[-R, R]$,
while Proposition \ref{prop76} provides a limit curve $c:[-R, R] \to Y$ such that $f\circ c= c_v |[-R, R]$ and $d(y, c(t)) = |t|$ for all $t\in [-R, R]$.
This implies that $c(r) \in \partial B(y, r)$ and $f\circ c(r)= c_v (r)$, hence $f(\partial B(y, r))\supseteq S_P(r)$. Similarly, we obtain $f(B(y, r))\supseteq N_P(r)$ for all $r>0$.
\kasten\vspace{0.2cm}

The following lemma is a consequence of Proposition \ref{prop63}. In this section we let $k>0$ denote an upper bound for the absolute values of the sectional curvatures of $M$.
\begin{lemma}\label{lemma135}
Let $y\in Y$ and $r\in (0, r_0)$, where $r_0>0$ is given by Proposition \ref{prop63}. Then there exists a closed, 1-Lipschitz curve $\gamma:[0, l]\to Y$ such that
$\gamma([0, l]) \supseteq \partial B(y, r)$ and $l\leq 2\pi\frac{\sinh(\sqrt{k}r)}{\sqrt{k}}$.
\end{lemma}

\proof
We choose a sequence $q_i\in F_i$ converging to $y$. According to \cite{Har64}, Lemma 5.2, we can choose a sequence $r_i$ of non-exceptional values of the distance functions
$d_i(q_i, \cdot)$ such that $\lim_{i\to \infty} r_i = r$. Then \cite{Har64}, Proposition 6.1, see also \cite{SST03}, Theorem 4.4.1, together with Proposition \ref{prop63} imply that
$\partial B(q_i, r_i)$ is a piecewise smooth, simple, closed curve. We let $\gamma_i: [0, l_i]\to \partial B(q_i, r_i)$ be a parametrization of $\partial B(q_i, r_i)$ by arclength,
in particular $l_i = {\cal H}^1(\partial B(q_i, r_i))$. Then (\ref{eq131}) and (\ref{eq113}) imply
$$
\limsup_{i\to \infty}l_i \leq 2\pi\frac{\sinh(\sqrt{k}r)}{\sqrt{k}}.
$$
Choosing a subsequence we may assume that the sequence $l_i$ converges to $l\leq 2\pi\frac{\sinh(\sqrt{k}r)}{\sqrt{k}}$ and, by the Arzel\`a-Ascoli theorem, that the $\gamma_i$
converge uniformly to a closed, 1-Lipschitz curve $\gamma:[0, l] \to Y$. It remains to be shown that $\partial B(y, r) \subseteq \gamma([0, l])$.
To prove this we argue by contradiction, and assume:
\begin{align}
\label{eq135} \text{There exists } z\in \partial B(y, r) \text{ and } \delta>0 \text{ such that } \delta_i(\partial B(q_i, r_i),z)\geq \delta\\
\nonumber \text{ for infinitely many } i\in\N.
\end{align}

Since $z\in \partial B(y, r)$ there exists $z'\in Y$ such that $d(y, z')> r$ and $d(z, z')<\delta/2$, in particular $d(y, z') < r+\delta/2$.
Now we choose a sequence $x_i'\in F_i$ converging to $z'$, and estimate
$$
|d_i(q_i,x_i')-d(y, z')| \leq \delta_i(q_i, y)+\delta_i(z', x_i'),
$$
where $r+\delta/2 > d(y, z') >r = \lim_{i\to \infty} r_i$, and $\lim_{i\to \infty} \delta_i(q_i, y)=0=\lim_{i\to \infty} \delta_i(z', x_i')$. Hence, for almost all $i\in \N$, we have
$$
r_i < d_i(q_i, x_i') < r + \delta/2.
$$
For these $i\in \N$ we choose a shortest geodesic from $q_i$ to $x_i'$ and, on this geodesic, the point $x_i$ with $d_i(q_i, x_i)=r_i$. Then $x_i\in \partial B(q_i, r_i)$, and
$$
\delta_i(x_i, z) \leq d_i(x_i, x_i') + \delta_i(x_i', z')+d(z', z).
$$

Since $d_i(x_i, x_i') = d_i(q_i, x_i') - r_i < r- r_i+ \delta/2$, $\lim_{i\to \infty} \delta_i(x_i', z')=0$, and $d(z', z) < \delta/2$, we obtain $\delta_i(x_i, z)<\delta$ for almost all
$i\in \N$, in contradiction to (\ref{eq135}). So, we have $\lim_{i\to \infty} \delta_i(\partial B(q_i, r_i), z)=0$ for all $z\in \partial B(y, r)$. Hence, for every $z\in \partial B(y, r)$,
we can find a sequence $s_i\in [0, l_i)$ such that $\lim_{i\to \infty} \gamma_i(s_i)=z$. Then a subsequence of the $s_i$ converges to some $s\in [0, l]$ and, by the uniform convergence of $\gamma_i$ to $\gamma$, we have $\gamma(s)=\lim_{i\to \infty} \gamma_i(s_i)=z$.
\kasten\vspace{0.2cm}

Lemma \ref{lemma135} has the following important consequences.

\begin{lemma}\label{lemma136}
Suppose $q_i\in F_i$ is a good sequence converging to $y\in Y$. Then the sequence $df_i(T_{q_i} F_i)$ converges in $G_2M$.
\end{lemma}

\proof Since $G_2M$ is compact it suffices to show that the sequence $df_i(T_{q_i} F_i)$ cannot have two different limit planes $P\neq P'$. If this were the case we would have
$f(\partial B(y, r))\supseteq S_P(r) \cup S_{P'}(r)$ for $r\in (0, \injrad(M))$, see Lemma \ref{lemma134}. Since $f$ is 1-Lipschitz and $P\neq P'$ this would imply
$$
\liminf_{r\to 0} \frac{{\cal H}^1(\partial B(y, r))}{2\pi r}\geq 2,
$$
while Lemma \ref{lemma135} implies
$$
\limsup_{r\to 0} \frac{{\cal H}^1(\partial B(y, r))}{2\pi r}\leq 1.
$$
\kasten\vspace{0.2cm}

According to (\ref{eq134}) and Lemma \ref{lemma136} the following definition makes sense.

\begin{definition}\label{def137}
We define $\tilde{f}:(Y, y_0)\to (G_2M, P_0)$ by $\tilde{f}(y)=P$ iff $\lim_{i\to \infty} df_i(T_{q_i} F_i)=P$ for every good sequence $q_i\in F_i$ with
$\lim_{i\to \infty} q_i=y$. In particular, $\tilde{f}$ is a lift of $f$, i.\,e. $\pi_G\circ\tilde{f}=f$.
\end{definition}

To formulate the next lemma we choose $r_1>0$ such that
\begin{equation}\label{eq136}
2\pi\frac{\sinh(\sqrt{k}r)}{\sqrt{k}} - 2\pi\frac{\sin(\sqrt{k}r)}{\sqrt{k}}< 2r\quad\mbox{for } r\in(0, r_1),
\end{equation}
and set
\begin{equation}\label{eq137}
r_2 = \min\Big\{r_0, r_1, \injrad(M), \frac{\pi}{\sqrt{k}}\Big\}.
\end{equation}

\begin{lemma}\label{lemma138}
If $y, z\in Y$ and $0<d(y, z)<r_2$, then $f(y)\neq f(z)$.
\end{lemma}

\begin{remark}\label{rem139}
{\em Once Proposition \ref{prop133} will be proven, we will know that $f(y)\neq f(z)$ under the weaker condition $0<d(y, z)<2\,\injrad(M)$. }
\end{remark}

\proof We argue by contradiction and assume that $y, z\in Y$, $0<d(y, z)< r_2$, but $f(y)=f(z)$. We set $d(y, z) = r$ and, in a first step, we additionally assume that $z\in\partial B(y, r)$. Lemma \ref{lemma135} provides a closed, 1-Lipschitz curve $\gamma: [0, l]\to Y$ such that $\gamma([0, l])\supseteq \partial B(y, r)$ and
$l\leq 2\pi\frac{\sinh(\sqrt{k}r)}{\sqrt{k}}$. From Lemma \ref{lemma134} we conclude that $(f\circ \gamma)([0, l])\supseteq S_P(r)$, where $P=\tilde{f}(y)$. Since $f(z) = f(y)$ we know that $d^M(f(z), S_P(r)) = r$. Using our assumption $z\in \partial B(y, r)$, we obtain
\begin{equation}\label{eq138}
2\pi\frac{\sinh(\sqrt{k}r)}{\sqrt{k}} \geq l \geq {\mathrm{length}}(f\circ \gamma)\geq {\cal H}^1(S_P(r)) + 2r.
\end{equation}
Since $r<\min\{\injrad(M), \frac{\pi}{\sqrt{k}}\}$, the Rauch comparison theorem implies
\begin{equation}\label{eq139}
{\cal H}^1(S_P(r))\geq 2\pi\frac{\sin(\sqrt{k}r)}{\sqrt{k}}.
\end{equation}
Since $0 < r < r_1$, inequalities (\ref{eq138}) and (\ref{eq139}) contradict (\ref{eq136}). Finally we treat the case that $f(y)=f(z)$, $0< d(y,z)=r< r_2$, but
$z\notin \partial B(y, r)$. Since $Y$ is a length space we can find a shortest connection in $Y$ from $y$ to $z$, and on this shortest connection, a sequence of points $z_n \neq z$
converging to $z$. Then $z_n \in \partial B(y, r_n)$, where $r_n = d(y, z_n)$, $r_n \ua r$, and $d^M(f(z_n), f(z)) \leq r-r_n$. Now we repeat the argument used above with
$\partial B(y, r)$ replaced by $\partial B(y, r_n)$, and note that
$$
d^M(f(z_n), S_P(r_n)) \geq r_n-(r-r_n).
$$
So we obtain
$$
2\pi\frac{\sinh(\sqrt{k}r)}{\sqrt{k}} \geq  2\pi\frac{\sin(\sqrt{k}r)}{\sqrt{k}} + 2(2r_n - r)
$$
and, for $n\to \infty$, the same contradiction as before.
\kasten\vspace{0.2cm}

Our next aim is to prove that, for sufficiently small $r>0$ and for all $y\in Y$, $P=\tilde{f}(y)$, we have
$$
f(B(y, r)) = N_P(r) = \{c_v(t)|(v, t)\in SP\times[0, r)\},
$$
see Corollary \ref{coroll1312}. We set $r_c=\min\{\convrad(M), \frac{1}{2}\injrad(M)\}$, where $\convrad(M)$ denotes the convexity radius of $M$. For every $p\in M$ and every $r\in(0, r_c)$ we have a smooth function $b_{p, r}$, defined on the unit tangent bundle over $B(p, r)$ by
$$
b_{p, r}(v) = \sup\{t>0| c_v(s)\in B(p, r)\mbox{ for all } s\in [0, t]\}\in (0, 2r).
$$
A {\em totally geodesic disk in} $B(p,r)$ is a 2-dimensional, totally geodesic, connected submanifold $N$ of $B(p, r)$ such that $\bar{N}\backslash N\subset \partial B(p, r)$.
So, if $N$ is a totally geodesic disk in $B(p, r)$ and $q\in N$ then
$$
N= \{c_v(t)|v\in S_q N, 0\leq t < b_{p,r}(v)\}.
$$
If $z\in f^{-1}(B(p,r))$ and $Q=\tilde{f}(z)$, then
$$
N(z, p, r) = \{c_v(t)|v\in SQ, 0\leq t < b_{p, r}(v)\}
$$
is the {\em totally geodesic disk in} $B(p, r)$ {\em determined by }$z$. We note the following facts concerning these notions.
\begin{align}
            &\text{If }N\neq\tilde{N} \text{ are totally geodesic disks in\,} B(p, r) \text{ and }q\in N\cap\tilde{N}, \text{ then }T_q N\neq T_q\tilde{N}. \label{eq1310} \\
\nonumber   &\text{Moreover, }{\cal H}^2(N\cap\tilde{N}) = 0.\\
            &\text{If }0<r<\tilde{r}<r_c \text{ and }z\in f^{-1}(B(p, r)), \text{ then }N(z, p, r) = N(z, p, \tilde{r})\cap B(p,r). \label{eq1311}
\end{align}
As a consequence of the Rauch comparison theorem, see (\ref{eq139}), we have:

\begin{align}
\label{eq1312} \text{If }P\in G_2M \text{ and }r\leq \min\{\injrad(M), \pi/\sqrt{k}\}, \text{ and }\\
\nonumber      N_P(r) =\{c_v(t)|(v, t)\in SP\times[0, r)\},\text{ then }\\
\nonumber      {\cal H}^2(N_P(r))\geq \frac{2\pi}{k}(1-\cos(\sqrt{k}r)) = a(r) >0.
\end{align}

From Proposition \ref{prop73} we know that ${\cal H}^2(B(y, r)) <\infty$ for every $y\in Y$, $r>0$. This implies the following preliminary result:

\begin{lemma}\label{lemma1310} Suppose $y\in Y$, $f(y) = p$, and $0<r<r_c$. Then there exist $n\in \N$ and $z_1,\ldots,z_n\in B(y, r)$ such that
$$
f(B(y, r))\subseteq \bigcup_{j=1}^n N(z_j, p, r).
$$
\end{lemma}

\proof We will show that $\{N(z, p, r)|z\in B(y, r)\}$ is a finite set. This will prove our claim since obviously $f(B(y, r))\subseteq \bigcup_{z\in B(y, r)} N(z, p, r)$.
We choose $\rho>0$ such that $r+\rho < r_c$, and assume that $w_1, \ldots, w_k\in B(y, r)$ are such that $N(w_1, p, r),\ldots, N(w_k, p, r)$ are pairwise different.
We set $\tilde{f}(w_1) = P_1, \ldots, \tilde{f}(w_k)= P_k$. From Lemma \ref{lemma134} we know that $N_{P_j}(\rho) \subseteq f(B(w_j, \rho))\subseteq f(B(y, r+\rho))$.
Moreover $N_{P_j}(\rho) \subseteq N(w_j, p, r+\rho)$. Using (\ref{eq1311}) we see that also the $N(w_j, p, r+\rho)$, $1\leq j\leq k$, are pairwise different.
Hence (\ref{eq1310}) implies ${\cal H}^2(N_{P_j}(\rho)\cap N_{P_{j'}}(\rho))=0$ if $1\leq j < j'\leq k$. Using (\ref{eq1312}) we obtain
$$
{\cal H}^2(B(y, r+\rho))\geq {\cal H}^2(f(B(y, r+\rho))\geq \sum_{j=1}^k {\cal H}^2(N_{P_j}(\rho)\geq k a(\rho).
$$
Hence we have $\# \{N(z, p, r)|z\in B(y, r)\} \leq a(\rho)^{-1} {\cal H}^2(B(y, r+\rho))<\infty$.
\kasten

\begin{lemma}\label{lemma1311}
Suppose $N$ is a totally geodesic surface in $M$, $V\subseteq Y$ is open and $f(\bar{V})\subseteq N$. Then $\tilde{f}(z) = T_{f(z)} N$ for all $z\in \bar{V}$.
\end{lemma}

\proof First note that Lemma \ref{lemma134} implies that $\tilde{f}(y) = T_{f(y)} N$ for all $y\in V$. If $z\in \bar{V}$ we can find $r\in(0, r_2/2)$ such that for every $q\in N$
with $d^M(q, f(z))< r$ there exists $v\in S_qN$ such that $c_v(t) = f(z)$ where $t=d^M(q, f(z))$. Since $z\in \bar{V}$ there exists $y\in V\cap B(z, r)$. Then $f(y)\in N$ and
$d^M(f(y), f(z)) < r$, so that we can find $v\in S_{f(y)} N$ such that $c_v(t) = f(z)$ if $t = d^M(f(y), f(z))$. Now let $q_i\in F_i$ be a very good sequence converging to $y$.
Then there exists a sequence $(v_i, t_i)\in S_{q_i} F_i\times \R$ such that $c_{v_i}(t_i)$ is a good sequence, and such that $\lim_{i\to \infty}(df_i(v_i), t_i) = (v, t)$ and
$\lim_{i\to \infty}\int_{-t-1}^{t+1}|A_i|^2\circ c_{v_i}(s)\,ds=0$. Then Proposition \ref{prop23} implies that $\lim_{i\to \infty}f_i\circ c_{v_i}(t_i) = c_v(t) = f(z)$, and that
\begin{equation}\label{eq1313}
\lim_{i\to \infty}df_i(T_{c_{v_i} (t_i)}F_i) = {\cal P}_t^{c_v}(\tilde{f}(y))= {\cal P}_t^{c_v}(T_y N) = T_{f(z)}N.
\end{equation}
Choosing a subsequence we may assume that the sequence $c_{v_i}(t_i)$ converges to some $w\in Y$. Since $d(y, w) \leq \lim_{i\to \infty} t_i = t$, we have $d(z, w) < r+t < r_2$.
Moreover $f(w) = \lim_{i\to \infty}f_i(c_{v_i}(t_i)) = f(z)$, so that $w=z$ by Lemma \ref{lemma138}. Hence (\ref{eq1313}) implies $\tilde{f}(z) = \tilde{f}(w) = T_{f(z)} N$.
\kasten

\begin{corollary}\label{coroll1312}
Suppose $y\in Y$, $P=\tilde{f}(y)$, and $0<r<r_c$. Then $f(B(y, r))=N_P(r)$.
\end{corollary}

\proof From Lemma \ref{lemma1310} we obtain totally geodesic disks $N_1, \ldots, N_n$ in $B(f(y), r)$ such that $f(B(y,r))\subseteq \bigcup_{j=1}^n N_j$. Obviously, we may assume that
$N_j\neq N_k$ for $1\leq j<k\leq n$. Since the $N_j$ are closed subsets of $B(f(y), r)$, the sets $V_j = \{z\in B(y, r)|f(z)\in N_j\backslash \bigcup_{k\neq j} N_k\}$ are open in $B(y,r)$,
and $f(B(y, r)\backslash \bigcup_{j=1}^n V_j)\subseteq \bigcup_{j<k} (N_j \cap N_k)$. Since ${\cal H}^2(N_j\cap N_k)=0$ by (\ref{eq1310}), and ${\cal H}^2(f(W))>0$ for every open set
$\emptyset \neq W\subseteq Y$ by Lemma \ref{lemma134}, we see that $B(y, r)\subseteq \bigcup_{j=1}^n \bar{V_j}$. Since $Y$ is a length space, the metric ball $B(y, r)$ is connected.
Hence, if $n>1$ there exist $1\leq j<k\leq n$ and $z\in \bar{V_j}\cap \bar{V_k} \cap B(y, r)$. Now Lemma \ref{lemma1311} implies that $\tilde{f}(z) = T_{f(z)} N_j = T_{f(z)} N_k$,
in contradiction to $N_j\neq N_k$ and (\ref{eq1310}). Hence we have $n=1$ and $f(B(y, r))\subseteq N_1$. Since $f(y) \in N_1$ and $P= \tilde{f}(y)$, we see that $N_1= N(P, f(y), r) = N_P(r)$. Finally, Lemma \ref{lemma134} implies that $N_P(r) \subseteq f(B(y, r))$, so that indeed $f(B(y, r))=N_P(r)$.
\kasten

\begin{remark}\label{rem1313}{\em It is conceivable that the proof of Corollary \ref{coroll1312} can be shortened by use of results from the theory of spaces of bounded integral curvature.
In particular, the estimate in \cite{Deb20}, Corollary 3.2(1), would easily imply Corollary \ref{coroll1312}. However, due to the global assumptions in \cite{Deb20},
it is not directly applicable in our situation.}
\end{remark}

\textsc{Proof of Proposition \ref{prop133}:} We will prove that the lift $\tilde{f}:(Y, y_0)\to (G_2M, P_0)$ of $f$ given by Definition \ref{def137} is a locally isometric covering map onto the leaf $L$ of ${\cal D}_2(M, g)$ through $P_0$. First note that, by Lemma \ref{lemma1311} and Corollary \ref{coroll1312}, the sets $\tilde{f}^{-1}(L)$ and
$\tilde{f}^{-1}(G_2 M\backslash L)$ are both open in $Y$. Since $Y$ is connected and $y_0\in \tilde{f}^{-1}(L)$ we see that $\tilde{f}(Y)\subseteq L$.
Let $d^L$ denote the distance on $L$ induced by $(\pi_G|L)^{\ast}g$. We will prove that $d(y, z) = d^M(f(y), f(z)) = d^L(\tilde{f}(y), \tilde{f}(z))$, if $y, z\in Y$ and
$d(y, z) < \min\{r_2/2, r_c\}$. We set $P=\tilde{f}(y)$ and $\rho= d^M(f(y), f(z))$. Then Corollary \ref{coroll1312} implies that $f(z)\in S_P(\rho)$. From Lemma \ref{lemma134} we obtain $w\in\partial B(y, \rho)$ such that $f(w) = f(z)$. Since $d(z, w) \leq d(z, y)+\rho\leq 2d(y, z) < r_2$, Lemma \ref{lemma138} implies $z=w$, hence $d(y, z)=d(y, w)=\rho= d^M(f(y), f(z))$.
Using this and Lemma \ref{lemma134} we see that, for $0<r<\frac{1}{2}\min\{r_2/2, r_c\}$, and for every $y\in Y$, $f|B(y, r)$ is an isometry from $B(y, r)$ onto $N_P(r)$ where $P= \tilde{f}(y)$.
Finally, since $r< r_c \leq \injrad(M)$ we know that $\pi_G|L_P(r)$ is an isometry from the metric ball $L_P(r)$ with center $P$ and radius $r$ in $L$ onto $N_P(r)$. Hence
$\tilde{f}|B(y, r) = (\pi_G|L_P(r))^{-1}\circ (f|B(y, r))$ is an isometry from $B(y, r)$ onto $L_P(r)$ with the distance $d^L$. Since $L$ is connected this implies that $\tilde{f}: Y\to L$ is a covering map.
\kasten

\textsc{Proof of Theorem \ref{thrm131}:} In view of Lemma \ref{lemma134} it suffices to derive a contradiction from the assumption that there exist $R>0$ and a sequence
$q_i\in \overline{B(p_i, R)}$ such that the sequence $f_i(q_i)$ has a limit point in $M\backslash \overline{N_{P_0}(R)}$. If this is the case we can choose a subsequence, denoted by the same symbols, such that all of the following statements hold:
\begin{itemize}
\item[(i)]    $\lim_{i\to \infty} f_i(q_i)$ exists and $\lim_{i\to \infty} f_i(q_i)\notin \overline{N_{P_0}(R)}$.
\item[(ii)]   The sequence $(F_i, p_i)_{i\in\N}$ converges to $(Y, y_0)$ with respect to (definite) pointed Gromov-Hausdorff convergence, and the sequence $(f_i)_{i\in \N}$
              converges to $f:Y\to M$.
\item[(iii)]  The sequence $(q_i)_{i\in\N}$ converges to some $y\in \overline{B(y_0, R)}$.
\end{itemize}
Then we have $f(y) = \lim_{i\to \infty} f_i(q_i)\notin \overline{N_{P_0}(R)}$. On the other hand, $(Y, d)$ is a length space and $\tilde{f}:Y\to L$ is locally isometric, so that
$\tilde{f}(\overline{B(y_0, R)})\subseteq \overline{L_{P_0}(R)}$. Hence $f(y) =\pi_G(\tilde{f}(y))\in \pi_G(\overline{L_{P_0}(R)})=\overline{N_{P_0}(R)}$, cf. Remark \ref{rem132}.
\kasten

We mention the following consequence of Proposition \ref{prop133}:
\begin{corollary}\label{coroll1314} Under the assumptions (\ref{eq131})--(\ref{eq133}) we consider the sets
$\Omega_i(R) = \{(x, x')|\{x, x'\}\subseteq B(p_i, R), d_i(x, x')\leq \injrad(M)\}$. Then
$$
\lim_{i\to \infty}\Big(\sup_{(x, x')\in\Omega_i(R)}\big|d_i(x, x')-d^M(f_i(x), f_i(x'))\big|\Big)=0
$$
for every $R>0$.
\end{corollary}

\proof  Otherwise we can find $R>0$ and a sequence $(x_i, x_i')\in\Omega_i(R)$ such that $\limsup_{i\to \infty}|d_i(x_i, x_i')-d^M(f_i(x_i), f_i(x_i'))|>0$.
Choosing a subsequence we may assume that actually $\lim_{i\to \infty}|d_i(x_i, x_i')-d^M(f_i(x_i), f_i(x_i'))|>0$ and, additionally, that the $f_i: F_i\to M$ converge to $f:Y\to M$,
and that $\lim_{i\to \infty} x_i = y\in Y$, $\lim_{i\to \infty} x_i' = y'\in Y$. Then we have $d(y, y')=\lim_{i\to \infty} d_i(x_i, x_i')\leq \injrad (M)$ and
$d^M(f(y), f(y'))=\lim_{i\to \infty} d^M(f_i(x_i), f_i(x_i'))$, and, hence $|d(y, y')- d^M(f(y), f(y'))|>0$. In contradiction to this inequality we will now show that
$d(y, y') = d^M(f(y), f(y'))$. Indeed, since $Y$ is a length-space there exists an arclength-parametrized curve $c:[0, d(y, y')]\to Y$ from $c(0)=y$ to $c(d(y, y'))= y'$.
Then Proposition \ref{prop133} implies that $\tilde{f}\circ c$ is a geodesic in $L$ and, hence, $f\circ c= \pi_G\circ \tilde{f}\circ c$ is a geodesic in $M$. This geodesic $f\circ c$
has length $d(y, y')\leq \injrad (M)$ and connects $f(y)$ to $f(y')$. This implies
$$
d(y, y')= {\mathrm{length}}(f\circ c)= d^M(f(y), f(y')).
$$
\kasten
Finally, we present a proof for Theorem \ref{thrm13}. It depends on Proposition \ref{prop133} and on an idea that is originally due to G. Reeb \cite{Reeb47}, namely Reeb's stability theorem.\\

\textsc{Proof of Theorem \ref{thrm13}:} We argue by contradiction and assume that there exists a sequence of complete, connected surface immersions $f_i:F_i\to M$ such that
$\volz^{i}(F_i)\to \infty$, while the sequence $E(f_i)$ is bounded. Under this assumption the well-known Vitali covering argument provides sequences $\tilde{p}_i\in F_i$ and $R_i\to \infty$
such that $\lim_{i\to \infty} E(f_i|B(\tilde{p}_i, R_i))=0$. Using Lemma \ref{lemma83} we obtain sequences $p_i\in F_i$, $\ve_i\da 0$ and $\rho_i\to \infty$ such that
(\ref{eq131}) and (\ref{eq132}) hold. Choosing a subsequence we can assume that $\lim_{i\to \infty} df_i(T_{p_i}F_i)=P_0$ exists, i.\,e. also (\ref{eq133}) is satisfied. Now Proposition \ref{prop133} provides a locally isometric covering map $\tilde{f}:(Y, y_0) \to (L, P_0)$ from a Gromov-Hausdorff limit space $Y$ onto the leaf $L$ of ${\cal D}_2(M, g)$ through $P_0$.
We want to show that $L$ is not homeomorphic to $S^2$ or $\R P^2$. Otherwise $Y$ would be compact. Since the $F_i$ are connected this would imply that there exists $D>0$ such that
${\mathrm{diam}}(F_i)\leq D$ and, by Corollary \ref{coroll62}, $\volz^{i}(F_i) \leq \frac{1}{k}(2\pi + E(f_i))(\cosh(\sqrt{k}D)-1)$
% $\volz^{i}(F_i) \leq \frac{1}{\pm k}(2\pi + E(f_i))(\cosh(\sqrt{\pm k}D)-1)$
for infinitely many $i\in \N$, in contradiction to our assumption $\volz^{i}(F_i)\to \infty$. Hence $\pi_G|L$ is a complete, totally geodesic immersion into $M$, where $L$ is a connected surface different from $S^2$ or $\R P^2$. This contradicts our assumption on $M$.
\kasten

\begin{remark}\label{rem1315} {\em Similarly we see that, under the assumptions made above, the lamination structure ${\cal L}$ on $S=\bar{L}$ from Proposition \ref{prop106} does not have
a leaf homeomorphic to $S^2$ or $\R P^2$.}
\end{remark}

%%%%%%%%%%% evtl. in den ersten Zeilen von z. B. Prop. 7.3, Lemma 8.2, 8.3. einsetzen
%\begin{flushleft}
%\end{flushleft}
%%%%%%%%%%%%%%%%%%%%%%%%

%%%%%%%%%%%%%%%%%%%%%%%%%%%%%%%%%%%%%%%%%%%%%%%%%%%%%%%%%%%%%%%%%%%%%%%%%%%

\end{document}